\documentclass{LMCS}

\def\doi{8 (3:07) 2012}
\lmcsheading%
{\doi}
{1--38}
{}
{}
{Mar.~\phantom05, 2012}
{Aug.~10, 2012}
{}

\usepackage{proof}
\usepackage{enumerate}
\usepackage{hyperref}
\usepackage{amsmath, amssymb}

\newcommand{\AND}[0]{\wedge}
\newcommand{\OR}[0]{\vee}
\newcommand{\IMPLY}[0]{\supset}
\newcommand{\ENTAIL}[0]{\rightarrow}

\DeclareMathOperator{\PIND}{PIND}
\DeclareMathOperator{\Env}{Env}
\DeclareMathOperator{\BdEnv}{BdEnv}
\DeclareMathOperator{\LEFT}{L}
\DeclareMathOperator{\RIGHT}{R}
\DeclareMathOperator{\LEFTb}{Lb}
\DeclareMathOperator{\RIGHTb}{Rb}

\DeclareMathOperator{\Prf}{\textup{-Prf}}
\DeclareMathOperator{\Con}{\textup{-Con}}
\DeclareMathOperator{\RCon}{\textup{RCon}}

\DeclareMathOperator{\PV}{\textup{PV}}
\DeclareMathOperator{\Cond}{\textup{Cond}}

\begin{document}

\title{Bounded arithmetic in free logic} 
\author[Y.~Yamagata]{Yoriyuki Yamagata}
\address{National Institute of Advanced
Industrial Science and Technology (AIST) \newline 3-11-46 Nakoji \\ Amagasaki, 661-0974 Japan}
\email{yoriyuki.yamagata@aist.go.jp}


\begin{abstract}
  One of the central open questions in bounded arithmetic is whether
  Buss' hierarchy of theories of bounded arithmetic collapses or not.
  In this paper, we reformulate Buss' theories using \emph{free
    logic} and conjecture that such theories are easier to handle.  To
  show this, we first prove that Buss' theories prove consistencies
  of induction-free fragments of our theories whose formulae have
  bounded complexity.  Next, we prove that although our theories are
  based on an apparently weaker logic, we can interpret theories in
  Buss' hierarchy by our theories using a simple translation.
  Finally, we investigate finitistic G\"{o}del sentences in our
  systems in the hope of proving that a theory in a lower level of
  Buss' hierarchy cannot prove consistency of induction-free
  fragments of our theories whose formulae have higher complexity.
\end{abstract}

\maketitle

\section{Introduction}

One of the central open questions in bounded arithmetic is whether
Buss' hierarchy $S^1_2 \subseteq T^1_2 \subseteq S^2_2 \subseteq T^2_2
\subseteq \cdots$ of theories of bounded arithmetic collapses
\cite{Buss:Book} or not.   Since it
is known that collapse of Buss' hierarchy implies the collapse of the
polynomial-time hierarchy \cite{Krajicek1991},
demonstration of the non-collapse of the theories in Buss' hierarchy could
be one way to establish the non-collapse of the polynomial-time hierarchy.
A natural way to demonstrate non-collapse of the theories in Buss'
hierarchy would be to identify one of these theories that proves (some
appropriate formulation of) a statement of the consistency of some
theory below it in the hierarchy.

Here, it is clear that we need a delicate notion of consistency
because of several negative results that have already been
established.  The ``plain'' consistency statement cannot be used to
separate the theories in Buss' hierarchy, since Paris and Wilkie
\cite{Wilkie1987Scheme} show that $S_2\; (\equiv \bigcup S^i_2)$
cannot prove the consistency of Robinson Arithmetic $Q$.  Apparently,
this result stems more from the use of predicate logic than from the
strength of the base theory.  However, Pudl\'ak \cite{PudlakNote}
shows that $S_2$ cannot prove the consistency of proofs that are
carried out within $S^1_2$ and are comprised entirely of bounded
formulae. Even if we restrict our attention to the induction-free
fragment of bounded arithmetic, we cannot prove the consistency of
such proofs, as shown by Buss and Ignjatovi\'c
\cite{BussUnprovability}.  More precisely, Buss and Ignjatovi\'c prove
that $S^i_2$ cannot prove the consistency of proofs that are comprised
entirely of $\Sigma^b_i$ and $\Pi^b_i$ formulae and use only BASIC
axioms (the axioms in Buss' hierarchy other than induction) and the
rules of inference of predicate logic.

Therefore, if we want to demonstrate non-collapse of the theories in
Buss' hierarchy, we should consider a weaker notion of consistency
and/or a weaker theory.  A number of attempts of this type have been
made, both on the positive side (those that establish provability of
consistency of some kind) and on the negative side (those that
establish non-provability of consistency).  On the positive side,
Kraj\'{\i}\v{c}ek and Takeuti \cite{Krajicek1992} show that $T^i_2
\vdash \RCon(T^i_1)$, where $T^i_1$ is obtained from $T^i_2$ by
eliminating the function symbol \#, and $\RCon(T^i_1)$ is a sentence
which states that all ``regular'' proofs carried out within $T^i_1$
are consistent.  Takeuti \cite{Takeuti1996}, \cite{Takeuti2000} shows
that there is no ``small'' strictly $i$-normal proof $w$ of
contradiction.  Here, ``$w$ is small'' means that $w$ has its
exponentiation $2^w$.  Although Takeuti allows induction in strictly
$i$-normal proof $w$, the assumption that $w$ is small is a
significant restriction to $w$ since bounded arithmetics cannot prove
existence of exponentiation.  Another direction is to consider
cut-free provability.  Paris Wilkie \cite{Wilkie1987Scheme} mentioned
above proves that $I\Delta_0 + \textup{exp}$ proves the consistency of
cut-free proofs of $I\Delta_0$.  For weaker theories than $I\Delta_0 +
\textup{exp}$, we need to relativized the consistency by some cut,
then we get similar results \cite{Pudlak1985}, \cite{Adamowicz2011}.
For further weaker theories, Beckmann \cite{Beckmann2002} shows that
$S^1_2$ proves the consistency of $S^{-\infty}_2$, where
$S^{-\infty}_2$ is the equational theory which is formalized by
recursive definitions of the standard interpretations of the function
symbols of $S_2$.  Also, it is known that $S^i_2$ proves
$\textup{Con}(G_i)$, that is, the consistency of quantified
propositional logic $G_i$.  On the negative side, we have the results
mentioned above, that is, those of Paris and Wilkie
\cite{Wilkie1987Scheme}, Pudl\'ak \cite{PudlakNote}, and Buss and
Ignjatovi\'c \cite{BussUnprovability}.  In addition, there are results
which extend incompleteness theorem to Herbrand notion of consistency
\cite{Adamowicz2001},\cite{Adamowicz2002}.

In this paper, we introduce the theory $S^i_2E$ ($i=-1,0,1,2\ldots$),
which for $i\ge 1$ corresponds to Buss' $S^i_2$, and we show that the
consistency of strictly $i$-normal proofs that are carried out only in
$S^{-1}_2E$, can be proved in $S^{i+2}_2$.  We improve on the
aforementioned positive results in that 1) unlike $T^i_1$ or $G_i$,
$S^i_2E$ is based on essentially the same language as $S^i_2$, thereby
making it possible to construct a G\"odel sentence by diagonalization;
2) unlike Takeuti \cite{Takeuti1996}, \cite{Takeuti2000}, we do not
assume that the G\"odel number of the proofs which are proved
consistent are small, that is, have exponentiations, thereby making it
possible to apply the second incompleteness theorem---in particular,
to derive a G\"odel sentence from the consistency statement; 3) unlike
the results on Herbrand and cut-free provability, $S^i_2E$ has the
Cut-rule, thereby, making it easy to apply the second incompleteness
theorem; and 4) unlike Beckmann \cite{Beckmann2002}, our system is
formalized in predicate logic.  On the other hand, we are still unable
to show that the consistency of strictly $i$-normal proofs is not
provable within $S^j_2$ for some $j\leq i$, but see Section
\ref{sec:goedel-sentence}.  In a sense, our result is an extension of
that of Beckmann \cite{Beckmann2002} to predicate logic, since both
results are based on the fact that the proofs contain ``computations''
of the terms that occur in them.  In fact, if we drop the Cut-rule
from $S^{-1}_2E$, the consistency of strictly $i$-normal proofs can be
proved in $S^1_2$ for any $i$.  This ``collapse'' occurs since,
roughly speaking, the combination of the Cut-rule and universal
correspond substitution rule in $\PV$.

$S^i_2E$ is based on the following observation: The difficulty in
proving the consistency of bounded arithmetic inside $S_2$ stems from
the fact that inside $S_2$ we cannot define the evaluation function
which, given an assignment of natural numbers to the variables, 
maps the terms of $S_2$ to their values. For example, the values of
the terms $2, 2\#2, 2\#2\#2,2\#2\#2\#2,\ldots$ increase exponentially;
therefore, we cannot define the function that maps these terms to
their values, since the rate of growth of every function which is
definable in $S_2$ is dominated by some polynomial in the length of
the input \cite{Parikh1971}.
With a leap of logic, we consider this fact to mean that we cannot
assume the existence of values of arbitrary terms in
bounded arithmetic. Therefore, we must explicitly prove the existence
of values of the terms that occur in any given proof.

Based on this observation, $S^i_2E$ is formulated by using \emph{free
  logic} instead of the ordinary predicate calculus.
Free logic is a logic which is free from ontological assumptions about
the existence of the values of terms.
Existence of such objects is explicitly stated by an
existential predicate rather than being
implicitly assumed.  See \cite{FreeLogic} for a general introduction
to free logic and \cite{Scott1979} for its application to
intuitionistic logic.

Using free logic, we can force each proof carried out within
$S^{-1}_2E$ to somehow ``contain'' the values of the terms that occur
in the proof. By extracting these values from the proof, we can
evaluate the terms and then determine the truth value of $\Sigma^b_i$
formulae.  The standard argument using a truth predicate proves the
consistency of strictly $i$-normal proofs that are carried out only in
$S^{-1}_2$.  It is easy to see that such a consistency proof can be
carried out in $S^{i+2}_2$.

The paper is organized as follows.  In Section \ref{sec:def} we define
$S^i_2E$\ and compare it to the systems of logic defined in
\cite{FreeLogic} and \cite{Scott1979}.  In Section
\ref{sec:consistency} we present the main result of this paper: a
proof inside $S^{i+2}_2$ of the consistency of induction-free strictly
$i$-normal proofs.  In Section \ref{sec:bootstrapping} we prove the
``bootstrapping theorem,'' in which we show that although our theory
$S^i_2E$ is based on an apparently weaker logic, we can interpret
$S^i_2$ inside $S^i_2E$ if $i\ge 1$.  Unlike the interpretation of
$S^i_2$ by $Q$ \cite{Hajek1998}, our interpretation does not increase
the complexity of formulae; in particular, all bounded formulae are
interpreted as bounded formulae.  Finally, in Section
\ref{sec:goedel-sentence} we raise the question of whether the
consistency of induction-free strictly $i$-normal proofs can be proved
inside $S^i_2$.  We consider a countably infinite set of finitistic
G\"{o}del sentences of $S^{-1}_2E$ to investigate this question.

\section{Definition of $S^i_2E$}
\label{sec:def}

$S^{-1}_2E \subseteq S^0_2E \subseteq S^1_2E \subseteq S^2_2E
\subseteq \cdots$ is a hierarchy of theories resembling Buss'
hierarchy $S^{-1}_2 \subseteq S^0_2 \subseteq S^1_2 \subseteq S^2_2
\subseteq \cdots$. (For purposes of comparison with our system, we
include $S^{-1}_2$and $S^0_2$ in Buss' hierarchy, where $S^{-1}_2$ is
defined as the theory which consists of the formulae that can be
proved from the BASIC axioms via the rules of inference
of predicate logic, and $S^0_2$ is the set of formulae
that can be proved from the BASIC axioms via the rules of inference
of predicate logic together with induction on
quantifier-free formulae.)  Our system is equipped with a predicate
$E$ which signifies existence of the values of terms.  In this
section, we introduce the theories $S^i_2E$ ($i \geq -1$) and their
languages, and we prove their basic properties.

\begin{defi}
  The \emph{theory $S^i_2E$} consists of the formulae of $S^i_2E$ that
  can be proved from the union of a finite set of logical axioms and a
  finite set $\mathcal{A}$ of proper axioms (as defined in Section
\ref{subsec:axiom}) via the rules of inference of free logic (with
some modifications) and (for $i\geq 0$) the $\Sigma^b_i\text{-}\PIND$
rule; the latter is an induction principle for $\Sigma^b_i$ formulae
which is based on the binary representations of the nonnegative
integers.  We do not explicitly specify what $\mathcal{A}$ is.
Instead, we make it extensible, and we specify certain conditions that
$\mathcal{A}$ must satisfy. We do this so that the motivation behind
the axioms will be more transparent. We also allow the set
$\mathcal{F}$ of function symbols to be extensible.
$S^i_2E(\mathcal{F}, \mathcal{A})$ denotes the individual theory
obtained by the function symbols in $\mathcal{F}$ and the proper
axioms in $\mathcal{A}$.

  An \emph{$S^i_2E$ proof} ($i = -1, 0, 1,
  \dots$) is a formal deduction in which only the rules of inference
  of $S^i_2E$ are used. Thus in an
  $S^i_2E$ proof, induction is restricted to
  application of the $\Sigma^b_i\text{-}\PIND$ rule; in the $i = -1$
  case, induction is not allowed at all.

  The theory $S^i_2E$ is thus the set of formulae $A$ of $S^i_2E$ for
  which an $S^i_2E$ proof of the sequent
  $\ENTAIL A$ exists, which we denote by $S^i_2E \vdash A$. We call
  such formulae \emph{theorems of
    $S^i_2E$}. We also use the notation $\vdash$ for sequents: $S^i_2E
  \vdash \Gamma \ENTAIL \Delta$.
\end{defi}

In Subsection \ref{subsec:lang} we define the language of $S^i_2E$.
In Subsection \ref{subsec:axiom} we describe the conditions that must
be satisfied by $\mathcal{A}$, and in Subsection
\ref{subsec:deduction} we introduce the rules of inference
of $S^i_2E$.  In Subsection \ref{subsec:freelogic} we
compare our system to known systems of free logic.

\subsection{Language of $S^i_2E$}
\label{subsec:lang}
The vocabulary of $S^i_2E$ is obtained from that of $S^i_2$ by adding
the unary predicate symbol $E$ and replacing the set of function
symbols of $S_2$ with an arbitrary but finite set $\mathcal{F}$ of
function symbols which denote polynomial-time computable functions.
The formulae of $S^i_2E$ are built up from atomic formulae by use of
the propositional connectives $\neg, \AND, \OR$; the bounded
quantifiers $\forall x \leq t, \exists x \leq t$; and the unbounded
quantifiers $\forall x, \exists x$.  Implication ($\IMPLY$) is omitted
from the language, and negation ($\neg$) is applied only to equality
$=$ and inequality $\leq$.  These restrictions appear essential to
prove consistency.  If there is implication (or negation applied to
arbitrary formulae) in $S^{-1}_2E$, $S^{-1}_2E$ allows induction
speedup \cite{BussUnprovability} \cite{PudlakNote}, therefore
$S_2^{-1}E$ polynomially interprets $S^i_2E$, $i \geq 0$.  This allows
to prove $Ef(n)$ for any polynomial time $f$ in $S_2^{-1}E$ by a
proof whose length is bounded by some fixed polynomial of length of
binary representation of $n$.  However, this contradicts the statement
of soundness (Proposition \ref{prop:soundness})

Since the standard interpretation of all function symbols in
$\mathcal{F}$ are a polynomial-time computable functions, all the
function symbols of $S_2E$ are definable in $S^1_2$, and we assume
that Cobham's recursive definitions of polynomial-time computable
functions are attached to the corresponding function symbols.

We sometimes identify the function symbols of $S_2E$ with their
standard interpretations. The distinction between the two types of
entities will be clear from the context.

\begin{defi}
  A set $\mathcal{F}$ of function symbols (for polynomial-time
  computable functions) is \emph{well grounded} if it satisfies the
  following conditions.
  \begin{enumerate}[(1)]

  \item $\mathcal{F}$ contains the $n$-ary constant zero function
    $0^n$ (the $n$-ary constant function whose value is $0$); the
    $n$-ary projection function $\textup{proj}^n_l$\ (the $n$-ary
    function that outputs the $l$th element in a sequence of length
    $n$), $k=1,\ldots,n$; and the so-called binary successor functions
    $s_0$ and $s_1$, where $s_0$ (resp.\ $s_1$) is the unary function
    defined by $s_0(a):=2a$ (resp.\ $s_1(a):=2a+1$).  Note that the
    binary representation of $2a$ (resp.\ $2a+1$) is obtained by
    appending 0 (resp.\ 1) to the binary representation of $a$, whence
    the moniker \emph{binary successor function}. 
  \item If $f \in \mathcal{F}$ is defined from functions $g,
  h_1, \dots, h_n$ by
    composition, then $g, h_1, \dots, h_n \in \mathcal{F}$.
  \item If $f \in \mathcal{F}$ is defined from functions $g, h_1, h_2$
    by recursion, that is, $f$ is
    defined by the equations
    \begin{align}
      f(0, \vec{x}) &= g(\vec{x}) \\
      f(s_0 x, \vec{x}) &= h_0(x, \vec{x}, f(x, \vec{x})) \\
      f(s_1 x, \vec{x}) &= h_1(x, \vec{x}, f(x, \vec{x}))
    \end{align}
    then $g, h_1, h_2\in \mathcal{F}$.  Cobham's limit recursion on
    notations can be written as above, providing $\Cond, \# \in
    \mathcal{F}$.  Since $\Cond, \#$ can be defined by recursion as
    above, for any finite set of polynomial time functions $f_1,
    \cdots, f_n$, there is a well-grounded $\mathcal{F}$ such that
    $f_1, \cdots, f_n \in \mathcal{F}$.
  \end{enumerate}
\end{defi}

If $\mathcal{F}$ is well grounded (which is true of
every set $\mathcal{F}$ of function symbols we consider in this
paper), we can define the \emph{definition degree $d(f)$} of each $f
\in \mathcal{F}$.
\begin{defi}[Definition degree]
\label{defn:defdeg} Let $f \in \mathcal{F}$.
\begin{enumerate}[(1)]
\item If $f$ is $0^n$, $\textup{proj}^n_l$, $s_0$, $s_1$, then $d(f) := 0$.
\item If $f$ is defined from $g,h_1,\dots,h_n$ by composition, then
  \begin{equation}
    d(f) := 1 + \max\{d(g),d(h_1), \dots, d(h_n)\}
  \end{equation}
\item If $f$ is defined from $g,h_1,h_2$ by recursion, then
  \begin{equation}
    d(f) := 1 + \max\{d(g),d(h_1), d(h_2)\}
  \end{equation}
\end{enumerate}
\end{defi} Induction on the definition degree of $f$ is used in
Subsection~\ref{subsec:bootstrap_I} to prove the totality of $f$ (i.e.,
that $Ea_1,\dots,Ea_n\ENTAIL Ef(a_1,\dots,a_n)$, where
$a_1,\ldots,a_n$ are used here as meta-symbols for variables of
$S^i_2E$).  Now we define the vocabulary of
$S^i_2E(\mathcal{F}, \mathcal{A})$.
\begin{defi}[Vocabulary] The \emph{vocabulary of
    $S^i_2E(\mathcal{F}, \mathcal{A})$} consists of the following
  symbols.
  \begin{desCription}
  \item\noindent{\hskip-12 pt\bf Constant symbols:}\ The only constant symbol of
  $S^i_2$ is 0.
  \item\noindent{\hskip-12 pt\bf Variables:}\ The variables of $S^i_2E$ are $x_1, x_2, \dots$.
    We often use $a, b$ and $a_1,a_2, \dots$, $b_1,b_2,
    \ldots$, $x, y, x_1, x_2,\dots$ as meta-symbols for variables, and we
    often denote a finite sequence of variables by $\vec{a}$,
    $\vec{a'}$, or $\vec{b}$.
  \item\noindent{\hskip-12 pt\bf Function symbols:}\ The function symbols of $S^i_2E(\mathcal{F},
    \mathcal{A})$ are the symbols in the finite set $\mathcal{F}$.
    For all $i\ge 0$, we can interpret $S^i_2$ in
    $S^i_2E(\mathcal{F},\mathcal{A})$, provided that $\mathcal{F}$
    contains the function symbols of $S_2$ (those for the unary
    functions $S$, $\lfloor\frac{\cdot}{2}\rfloor$, $|\cdot|$ and the
    binary functions $\#$, $+$, $\cdot$), where $S$ is the successor
    function division by two rounded to 0, length defined by
    $S(a):=a+1$ and $\lfloor \frac{\cdot}{2} \rfloor$ is the function
    where $\left\lfloor \frac{a}{2}\right\rfloor$ is defined as the
    natural number $n$ such that $a\in\{2n,2n+1\}$. As stated earlier,
    $|a|$ is defined as the number of bits in the binary
    representation of $a$ (by convention, $|0|=0$). The binary
    function $\#$ is the so-called smash function defined by $a \# b
    := 2^{|a||b|}$, and $+$ and $\cdot$ are the usual addition and
    multiplication functions, respectively.  We assume that these
    function symbols are contained in $\mathcal F$.  Further as we
    exploit the binary representations of the natural numbers (as
    finite bit strings), we introduce the binary function $\oplus$ and
    $\ominus$.  $\oplus(a,b)$ is defined as the natural number whose
    binary representation is the concatenation of the binary
    representations of the natural numbers $a$ and $b$ (the bits of
    $a$ being the most significant bits of $\oplus(a,b)$).
    $\ominus(a, b)$ is defined as $\lfloor \frac{a}{2^{|b|}} \rfloor$.
    We also assume $\oplus$ and $\ominus$ are contained in $\mathcal
    F$.

    From this point on, we write them in infix notation (as
    in $u\oplus r$ and $u\ominus r$). Since the functions $\ominus$
    and $\oplus$ are polynomial-time functions, and thus
    $\Sigma^b_1$-definable in $S^1_2$, we use the notations $u\oplus
    r$ and $u\ominus r$ in our informal proofs as well as in actual
    formulae of $S_2$ and $S_2E$.
  \item\noindent{\hskip-12 pt\bf Predicate symbols:}\ $S^i_2E$ has three predicate symbols: $E$,
    $=$, $\leq$.  The unary predicate $E$ signifies that, for
    every term $t$ of which $E$ is asserted to hold, the value
    of $t$ actually exists(i.e., that it converges
    to a standard natural number). The binary predicate $=$
    signifies equality, and $\leq$ signifies the less-than-or-equal-to
    relation.  $p$ is used as a meta-variable for the predicate
    symbols $=$ and $\leq $.
  \end{desCription}
\end{defi}

\begin{defi}[Terms]
 The \emph{terms of $S^i_2E$} are defined recursively as follows.
 \begin{iteMize}{$\bullet$}
 \item $0$ and the variables $x_1,x_2,x_3,\dots$ are terms.
 \item If $f$ is an $n$-ary function symbol and $t_1, \dots, t_n$ are
   terms, then $f(t_1, \dots, t_n)$ is a term.
 \end{iteMize}
\end{defi} We use $s,t,t_1, t_2, \dots, u_1, u_2, \dots$ as
meta-variables for terms, and we often denote a finite sequence of
terms by $\vec{s}$ or $\vec{t}$.

For the functions $S, s_0, s_1$, we omit parentheses, instead denoting
$S(t)$, $s_0(t)$, and $s_1(t)$ by $St$, $s_0t$, and $s_1t$,
respectively.  Also, we write $\lfloor \frac{t}{2} \rfloor$ or $\lfloor
\frac{1}{2} t \rfloor$ for $\lfloor\frac{\cdot}{2}\rfloor t$, and we write
the binary functions $\#$, $+$, $\cdot$ in infix notation, as in $t +
u$.

We define \emph{numerals} as terms which are constructed from $s_0$,
$s_1$, and 0 alone.  They are used in $S^i_2E$ to denote natural
numbers.  We (informally) use the numbers $0, 1, 2, \dots$
(represented in the decimal system) to denote numerals.  In general,
we use the same notation for numerals as for the corresponding natural
numbers.  The distinction between the two types of entities will be
obvious from the context.

\begin{defi}[Formulae]\label{defn:formulae}
  The \emph{formulae of $S^i_2E$} are defined recursively as follows.
  \begin{iteMize}{$\bullet$}
\item If $t$ is a term, then $Et$ is a formula.
\item If $p$ is $=$ or $\leq$ and $t_1,t_2$ are terms, then
  $p(t_1,t_2)$ is a formula.
  \item If $p$ is $=$ or $\leq$ and $t_1,t_2$ are terms, then $\neg
    p(t_1, t_2)$ is a formula.
  \item If $A$ and $B$ are formulae, then $A \AND B$ and $A \OR B$ are
    formulae.
  \item If $A(a)$ is a formula and $x$ is a variable, then
    $\forall x. A(x)$ and $\exists x. A(x)$ are formulae.
  \item If $A(a)$ is a formula, $x$ is a variable, and $t$ is a term,
    then $\forall x \leq t. A(x)$ and $\exists x \leq t. A(x)$ are
    formulae.
  \end{iteMize}
\end{defi} We use $A, B, \dots, A_0,A_1, A_2, \dots$ as
meta-variables for formulae, and
$\Gamma,\Delta,\Pi,\Lambda,$ and $\Gamma_1,\Gamma_2,\dots,\Delta_1,\Delta_2,\dots$
as meta-variables for finite sequences of formulae. $\Gamma,\Pi$
denotes the concatenation of $\Gamma$ and $\Pi$ (in that order:
$\Gamma$ followed by $\Pi$). The concatenation of a finite sequence
$\Gamma$ and a single formula $A$ is denoted by $\Gamma,A$, while the
concatenation of $A$ and $\Gamma$ is denoted by $A,\Gamma$.

If $p$ is $=$ (resp.\ $\leq$), we write the formula $p(t_1,t_2)$ in
infix notation, as $t_1=t_2$ (resp.\ $t_1\leq t_2$).  We write the
negations of $t_1=t_2$ and $t_1\leq t_2$ as $t_1\not= t_2$ and
$t_1\not\leq t_2$, respectively. In the meta-language, a notation of the
form $a \equiv b$ means that $a$ and $b$ are the same syntactic
construct (the same variable, the same term, the same formula, or the
same sequence of variables/terms). We often write $E\vec{a}$ to denote
the sequence $Ea_1, \dots, Ea_k$ for a finite sequence $\vec{a}\equiv
a_1,\dots,a_k$ of variables.

Formulae have the usual meaning, and free and bound variables are
defined as usual.  We sometimes write $A(a)$ to indicate that $a$ is
possibly among the variables that occur free in $A$. This notation
does not imply that $a$ actually occurs in $A$; what it indicates is
that if $a$ does occur in $A$, then $a$ occurs free in $A$. We tend to
use meta-symbols such as $x,y,x_1,x_2,\dots$ for variables that are
captured by an outer quantifier, as in $\forall x.A(x)$ or $\exists
x.A(x)$, and meta-symbols such as $a,b,a_1,a_2,\dots,b_1,b_2,\dots$
otherwise. We write $A(t)$ for the formula which is obtained by
substituting the term $t$ for $a$ in $A(a)$ (or for the variable $x$
in $A(x)$).  We follow the convention of assuming that variables which
occur bound in $A(a)$ and also occur in $t$ are renamed in $A(a)$
before such a substitution is made, so that the variables in $t$ are
not ``accidentally'' bound in $A(t)$.

We say that a quantifier is \emph{bounded} if it is of the form
$\forall x \leq t$ or $\exists x \leq t$.  On the other hand, the
quantifiers $\forall x$ and $\exists x$, without the bound $\leq t$,
are called \emph{unbounded}.

A formula of the form $t_1=t_2$ or $t_1\leq t_2$ is called
\emph{atomic}, and a formula of the form $Et$,
$t_1=t_2$, $t_1\leq t_2$, $t_1\not= t_2$, or $t_1\not\leq t_2$ is
called \emph{basic}.  Also, a formula of the form $Et$ is called an
$E$-form, but in contrast to the usual convention for atomic formulae
in predicate logic, an $E$-form is defined not to be atomic.

In $S^i_2E$, implication is absent and negation is restricted to
atomic formulae.  We interpret the negation $\neg A$ for a non
$E$-form formula $A$ as the De Morgan dual $\overline{A}$, and the
implication $A \subset B$ as $\overline{A} \OR B$.  As shown in
Section \ref{subsec:bootstrap_III}, $\neg$ and $\IMPLY$ defined in
this way satisfy the usual rules of inference of predicate logic,
provided that all the functions represented by function symbols of
$S_2E$ are provably total (that is, for each $n$-ary function symbol
$f$, we can prove $Ea_1, \ldots, Ea_n \ENTAIL Ef(a_1, \ldots, a_n)$).
This fact is proved in Proposition \ref{prop:total}.

Since implication is absent, we cannot code bounded universal
quantification by unbounded universal quantification together with
implication.  Therefore, we need bounded quantifiers that are
constructs in their own right.

By formalizing the notion of ``usual meaning,'' we obtain an
interpretation of the formulae of $S^i_2E$ in $S_2$.
\begin{defi}[Standard
  Interpretation]\label{defn:interpretation} We interpret an $S^i_2E$
  formula $A$ in $S_2$ by first replacing every subformula of $A$
  which is of the form $Et$ with the $S_2$ formula $t=t$, and then
  replacing the function symbols in the resulting formula with their
  $S^1_2$ definitions.  Therefore, we obtain an interpretation of the
  $S^i_2E$ formulae in $S^i_2$, and for $i\in\{-1,0\}$ we obtain an
  interpretation of the $S^i_2E$ formulae in $S^1_2$.
\end{defi} In Subsection \ref{subsec:translation}, we present an
interpretation of the $S^i_2$ formulae in
$S^i_2E(\mathcal{F},\mathcal{A})$ for any $\mathcal{F}$ that contains
the function symbols of $S_2$.

We say that a formula of $S^i_2E$ is \emph{bounded} if all of its
quantifiers are bounded.  The bounded formulae of $S^i_2E$ are
classified into hierarchies $\Sigma^b_j, \Pi^b_j$ ($j \in
\mathbb N$).
\begin{defi}[$\Sigma^b_j, \Pi^b_j$]
  The classes $\Sigma^b_j, \Pi^b_j$ of bounded
  formulae of $S^i_2E$ are defined recursively as follows.
  \begin{iteMize}{$\bullet$}
  \item $\Sigma^b_0 = \Pi^b_0$ is the set of quantifier-free formulae
    (formulae without quantifiers).
  \item If $A, B \in \Sigma^b_j$, then $A \AND B, A \OR B \in
    \Sigma^b_j$.  Similarly, if $A, B \in \Pi^b_j$, then $A \AND B, A
    \OR B \in \Pi^b_j$.
  \item If $A \in \Sigma^b_j$, then $A \in \Pi^b_{j+1}$.  Similarly,
    if $A \in \Pi^b_j$, then $A \in \Sigma^b_{j+1}$.
  \item If $A(a) \in \Sigma^b_j$, then $\forall x \leq t. A(x) \in
    \Pi^b_{j+1}$.  Similarly, if $A(a) \in \Pi^b_j$, then $\exists x
    \leq t. A(x) \in \Sigma^b_{j+1}$.
  \end{iteMize}
\end{defi} These hierarchies are used to control the strength of
mathematical induction in $S_2E$, just as in Buss' system.

\subsection{Axioms of $S^i_2E$}
\label{subsec:axiom}

In this subsection, we first discuss the conditions which all axioms
of $S^i_2E$ must satisfy. Then we introduce the logical axioms, and
finally we impose certain conditions on the axioms in
$\mathcal{A}$ that will allow us to interpret $S^i_2$ inside
$S^i_2E$\ for $i\ge 1$.

For our proof of consistency of strictly $i$-normal proofs to work,
the axioms of $S^i_2E$ must satisfy \emph{the boundedness conditions}.

 \begin{defi}\label{defn:boundedness}
   A sequent $\Gamma\ENTAIL\Delta$ \emph{satisfies the boundedness
     conditions} if it has the following three
   properties, where $\vec{a}$ are the variables that occur free
     in $\Gamma\ENTAIL\Delta$.
     \begin{enumerate}[(1)]
     \item  All the formulae that occur in $\Gamma\ENTAIL\Delta$ are
     basic.
   \item Every variable in $\vec{a}$ occurs free in at least one
     formula in $\Gamma$.
   \item There is a constant $\alpha\in\mathbb{N}$ such that
     \begin{equation}
       S^1_2 \vdash \max\{\vec{t}_\Delta(\vec{a})\} \le \alpha
       \cdot\max\{\vec{t}_\Gamma(\vec{a})\},
     \end{equation} $\vec{t}_\Delta(\vec{a})$ are
   the subterms of the terms that occur in $\Delta$ and
   $\vec{t}_\Gamma(\vec{a})$ are the subterms of the terms that occur in
   $\Gamma$ (for convenience, $\max\emptyset$ is defined to be
   1).  Since  the function
     symbols of $S_2E$ are definable in $S^1_2$, we can regard the terms in
     $\vec{t}_\Gamma(\vec{a})$ and $\vec{t}_\Delta(\vec{a})$ as terms of
     $S^1_2$, hence  we can regard
     $\max\{\vec{t}_\Gamma(\vec{a})\} \le
     \alpha\cdot\max\{\vec{t}_\Delta(\vec{a})\}$ as an $S^1_2$ formula.
     \end{enumerate}
\end{defi}

For our proof of consistency of $i$-normal proofs, we need the
boundedness conditions to hold of the axioms.  The third condition
states that, under any given valuation of the variables that occur
free in $\Gamma\ENTAIL\Delta$, the values of the subterms of the terms
in $\Delta$ cannot exceed the values of the subterms of the terms in
$\Gamma$ by the constant factor $\alpha$, and that this fact can be
proved in $S^1_2$.  In general, the constant $\alpha$ varies from one
axiom to another. Since $\mathcal{A}$ is finite, however, there is a
single constant $\alpha$ that applies to all the axioms in
$\mathcal{A}$. For the logical axioms and the proper axioms required
to interpret $S^i_2E$ in $S^1_2$, we can take $\alpha$ to be $4$, so
from this point on we assume that $\alpha = 4$.  This property plays a
crucial role in our consistency proof.  The reason for imposing the
first and second boundedness conditions is that we want to avoid
complexities that would otherwise arise in satisfying the third
boundedness condition.

The logical axioms are divided into \emph{$E$-axioms} and
\emph{equality axioms}.  The $E$-axioms establish the basic properties
of the $E$ predicate, and the equality axioms establish the basic
properties of equality.

\begin{defi}[$E$-axioms]\label{defn:E-axioms}
\begin{align}
  Ef(a_1, \dots, a_n) &\ENTAIL Ea_j &j = 1, \dots, n;
  f\in\mathcal{F} \label{eq:Ax-E-f}\\ p(a_1, a_2)&\ENTAIL Ea_j &j=1,2;
  p\in\{=,\leq\} \label{eq:Ax-E-p}\\ \neg p(a_1, a_2)&\ENTAIL
  Ea_j&j=1,2; p\in\{=,\leq\} \label{eq:Ax-E-not-p}
\end{align}
\end{defi} Axiom \eqref{eq:Ax-E-f} states that if the value of
$f$ exists, then the arguments of $f$ also exist.  Axioms
\eqref{eq:Ax-E-p} and \eqref{eq:Ax-E-not-p} state that if $=$, $\leq$,
or the negation of either of them holds of $a_1, a_2$, then $a_1, a_2$
exist.  That is, in the terminology of Scott \cite{Scott1979}, we
assume that all functions and predicates are \emph{strict}.

\begin{defi}[Equality axioms]\label{defn:Equality-axioms}
\begin{align}
  Ea &\ENTAIL a=a \label{eq:Ax-eq-refl}\\ a=b &\ENTAIL
  b=a\label{eq:Ax-eq-sym} \\ a=b, b=c &\ENTAIL
  a=c \label{eq:Ax-eq-trans}\\ a=b &\ENTAIL s_j a = s_j b &j =
  0,1 \label{eq:Ax-eq-succ'}\\ Ef(\vec{a}), \vec{a} = \vec{b} &\ENTAIL
  f(\vec{a}) = f(\vec{b})&f\in\mathcal{F} \label{eq:Ax-eq-fun}\\
  a_1=b_1,a_2=b_2,p(a_1, a_2)&\ENTAIL p(b_1, b_2)&p\in\{=,\leq\}
  \label{eq:Ax-eq}
\end{align}
\end{defi} Axioms \eqref{eq:Ax-eq-sym}, \eqref{eq:Ax-eq-trans},
\eqref{eq:Ax-eq-succ'}, and \eqref{eq:Ax-eq} are
standard, but Axiom \eqref{eq:Ax-eq-refl} needs an
explanation.  For our soundness proof, we would like to have the
property that when a closed formula $\phi$ is proved,
the values of the terms which
occur in $\phi$ are bounded by the code
for the proof of $\phi$.  Therefore, we cannot use just $t = t$ as
a substitution instance of an equality axiom with a closed term $t$,
since we would have no control over
the value of $t$.  To deal with this problem, we add
$Ea$ to the antecedent of the usual reflexive law of equality, which ensures
that the proof of $t = t$ contains the proof of $Et$,
and hence that the code for the proof of $t=t$ exceeds the
value of $t$.

Axiom \eqref{eq:Ax-eq-succ'} states that the binary successor
functions $s_0, s_1$ preserve equality.  For the other function
symbols, we have only a restricted preservation of equality,
$Ef(\vec{a}), \vec{a} = \vec{b}\ENTAIL f(\vec{a}) = f(\vec{b})$, as an
axiom, and not full preservation of equality, $\vec{a} = \vec{b}
\ENTAIL f(\vec{a})= f(\vec{b})$.  This is because,the
case of the reflexive law of equality, for sequences $\vec{t},\vec{u}$
of closed terms we want to bound the values of $f(\vec{t})$ and
$f(\vec{u})$ by the code for the proof of the equality $f(\vec{t}) =
f(\vec{u})$.  Without the condition $Ef(\vec{a})$ in the antecedent of
Axiom \eqref{eq:Ax-eq-fun}, we cannot ensure that this property will
hold.  For example, if $t_1,t_2,u_1,u_2$ are closed terms and $t_1 =
u_1, t_2=u_2 \ENTAIL t_1 \# t_2 = u_1 \# u_2$ is a substitution
instance of the usual equality axiom (Axiom \eqref{eq:Ax-eq-fun} with
$f\equiv \#$, $\vec{a}\equiv (a_1,a_2)$, and $\vec{b}\equiv
(b_1,b_2)$, but without $Ef(\vec{a})$ in the antecedent), the values
of $t_1 \# t_2$ and $u_1 \# u_2$ are not necessarily bounded by the
code for the proof of the equality $t_1 \# t_2 = u_1 \# u_2$.
However, we can prove full preservation of equality for function
symbols as a theorem, without assuming $Ef(\vec{a})$ (see Subsection
\ref{sec:bootstrap_II}), since we can prove that all the function
symbols of $S_2E$ represent total functions, that is, that for all $f
\in \mathcal F$, $E\vec{a} \ENTAIL Ef(\vec{a})$.  Actually, using the
$\Sigma^b_2\text{-}\PIND$ rule, we can prove full preservation of
equality for function symbols by induction on the definition degree of
$f$, without using Axiom \eqref{eq:Ax-eq-fun} at all.  However, we
retain Axiom \eqref{eq:Ax-eq-fun} since we want to have the equality
axioms in $S^1_2E$.

Axiom \eqref{eq:Ax-eq} states that the truth value of the predicates
$=$ and $\leq$ is preserved by equality.

Next, we discuss the proper axioms (the axioms in $\mathcal{A}$).  In
$S^{i+2}_2$ we can prove the consistency of strictly $i$-normal
proofs---which are carried out within $S^{-1}_2E(\mathcal{F},
\mathcal{A})$---for any $\mathcal{A}$ all of whose elements satisfy
the boundedness conditions. To interpret $S^i_2$ inside
$S^i_2E(\mathcal{F}, \mathcal{A})$, we need to have the function
symbols of $S_2$ in $\mathcal{F}$, and to have the following data
axioms, separation axioms and auxiliary axioms in $\mathcal{A}$.

\begin{defi}[Data axioms]\label{defn:data-axiom}
  These are the \emph{data axioms} introduced by Leivant
  \cite{Leivant2001}.
\begin{align}
    &\ENTAIL E0 \label{eq:Ax-data-0}\\ Ea &\ENTAIL Es_ja &j =
    0,1 \label{eq:Ax-data-suc}
  \end{align}
\end{defi}

\begin{defi}[Separation axioms]
\begin{align}
  s_j a = s_j b &\ENTAIL a=b &j = 0, 1 \label{eq:Ax-sep-sisi}\\ Ea, Eb
  &\ENTAIL s_0 a \not= s_1b \label{eq:Ax-sep-s0s1}\\ Ea &\ENTAIL s_1 a
  \not= 0 \label{eq:Ax-sep-s1=0}\\ s_0a = 0 &\ENTAIL a =
  0 \label{eq:Ax-sep-s0=0}
\end{align}
\end{defi}

\begin{defi}[Defining axioms]\label{defn:defnAx}  It is
  clear that by Cobham's definition of polynomial-time computable
  functions, the defining equations for each $n$-ary function symbol
  $f \in \mathcal{F}$ can be written in the form $f(u(a_1), a_2,
  \dots, a_n) = t(a_1, a_2, \dots, a_n)$ where
  $u(a_1)$ is one of $0,a_1,s_0a_1,s_1a_1$.  For each defining equation in
this form, the following axiom is in $\mathcal{A}$.
\begin{equation}
    \label{eq:function_definition} Ea_1, Ea_2, \dots, Ea_n, Et(a_1,
    a_2, \dots, a_n) \ENTAIL f(u(a_1), a_2, \dots, a_n) = t(a_1, a_2,
    \dots, a_n)
  \end{equation} We call these axioms the \emph{defining axioms}.
\end{defi}

\begin{defi}[Auxiliary axioms]\label{defn:auxiliary}
  The \emph{auxiliary axioms} are those generated from the BASIC
  axioms of $S_2$ by the following procedure.
  
  First, render these axioms free of logical connectives (see
  \cite{Takeuti1996}).  For example, the first BASIC axiom, $b \leq a
  \IMPLY b \leq S a$, is changed to the sequent $b \leq a \ENTAIL b
  \leq S a$.

  Then for each term $t$ that occurs in the succedent, we add
  $Et$ to the antecedent.  The axiom in the example becomes $Eb, ESa,
  b \leq a \ENTAIL b \leq S a$.

  Finally, simplify the resulting axioms by removing the unnecessary assumptions---taking care, however, to ensure that the
  final version of each axiom satisfies the boundedness conditions
  (Definition \ref{defn:boundedness}).  For example, we simplify $Eb,
  ESa, b \leq a \ENTAIL b \leq S a$ to $ESa, b \leq a \ENTAIL b \leq S
  a$, since $Eb$ follows from $b \leq a$ by $E$-axiom
  \eqref{eq:Ax-E-p} and, after this simplification, this auxiliary
  axiom still satisfies the boundedness conditions.
\end{defi}

From the data axioms and separation axioms, we can prove the relations
$n = m$ and $n \not= m$ that hold between numerals $n, m$ in their
standard interpretations.  Moreover, we can prove totality of every $f
\in \mathcal{F}$ (i.e., that $Ea_1, \dots, Ea_n \ENTAIL Ef(a_1, \dots,
a_n)$) from the logical axioms, the defining axioms, and
$\Sigma^b_0\text{-}\PIND$ (see Subsection \ref{subsec:bootstrap_I}).
Also, we can use the auxiliary axioms, together with other axioms and
the rules of inference of $S^1_2E$, to prove the
interpretation of the BASIC axioms of $S_2$ in $S^1_2E$ which is
defined in Subsection \ref{subsec:translation} (see Subsection
\ref{subsec:bootstrap_III}).

From this point on, we assume that $\mathcal{F}$ contains all the
  function symbols of $S_2$ and that $\mathcal{A}$ contains all the
proper axioms outlined above.  All of these proper axioms, as well as
the logical axioms, satisfy the boundedness conditions.

\subsection{Rules of inference of $S^i_2E$}
\label{subsec:deduction}

We formulate $S^i_2E$ proofs using the sequent calculus. The rules of
inference of $S^i_2E$ consist of the rules for predicate logic (but
with a modified $\RIGHT\neg$-rule; see \eqref{eq:Rneg}) plus the
$\PIND$ rule applied to $\Sigma^b_i$ formulae, which is called the
$\Sigma^b_i\text{-}\PIND$ rule.  The $\Sigma^b_i\text{-}\PIND$ rule
derives $\Gamma, Et \ENTAIL A(t), \Delta$ from $\Gamma \ENTAIL \Delta,
A(0)$; $A(a), \Gamma \ENTAIL \Delta, A(s_0a)$; and $A(a), \Gamma
\ENTAIL \Delta, A(s_1a)$.

If we can derive a sequent $\Gamma \ENTAIL \Delta$ from premises
$\Gamma_1 \ENTAIL \Delta_1, \dots, \Gamma_n \ENTAIL \Delta_n$, we
write
\begin{equation}
    \infer[R,]{\Gamma \ENTAIL \Delta}{
      \Gamma_1 \ENTAIL \Delta_1 \quad \cdots \quad \Gamma_n \ENTAIL
    \Delta_n }
\end{equation} where $R$ is the name of the rule. The name is often
omitted if it is obvious. All the rules of inference which appear 
in this paper have at most three premises.

The rules of inference are divided into the \emph{identity rule}, the
\emph{axiom rule}, the \emph{structural rules}, the \emph{logical
  rules}, the \emph{Cut rule}, and the \emph{$\PIND$ rule}.

The \emph{identity rule} is used to express that a formula $A$ implies
itself.

\paragraph{Identity rule:}\begin{equation}\infer[\textup{Id}]{A
\ENTAIL A}{}\end{equation}

The \emph{axiom rule} is used to derive a sequent which is a
\emph{substitution instance} of an axiom, that is, a sequent of the
form $\Gamma(\vec{t}(\vec{b}))\ENTAIL\Delta(\vec{t}(\vec{b}))$ where
$\Gamma(\vec{a})\ENTAIL\Delta(\vec{a})$ is an axiom, $\vec{a}$ are the
variables that occur in $\Gamma\ENTAIL\Delta$,
$\vec{t}$ are the terms that are being substituted for the variables
in $\vec{a}$, and $\vec{b}$ are the variables that occur in $\vec{t}$.

\paragraph{Axiom rule:}\begin{equation}
  \infer[\textup{Ax}]{\Gamma(\vec{t}(\vec{b}))\ENTAIL\Delta(\vec{t}(\vec{b}))}{}
\end{equation}

The \emph{structural} rules are defined as usual.

\paragraph{Weakening rules:}
  \begin{equation} \infer{A, \Gamma \ENTAIL
      \Delta}{\Gamma \ENTAIL \Delta} \end{equation}
  \begin{equation} \infer{\Gamma \ENTAIL
      \Delta, A}{\Gamma \ENTAIL \Delta} \end{equation}

\paragraph{Contraction rules:}
  \begin{equation} \infer{A, \Gamma \ENTAIL
      \Delta}{A, A, \Gamma \ENTAIL \Delta} \end{equation}
  \begin{equation} \infer{\Gamma \ENTAIL
      \Delta, A}{\Gamma \ENTAIL \Delta, A, A} \end{equation}

\paragraph{Exchange rules:}
  \begin{equation} \infer{\Gamma, B, A, \Pi \ENTAIL
      \Delta}{\Gamma, A, B, \Pi \ENTAIL \Delta} \end{equation}
  \begin{equation} \infer{\Gamma \ENTAIL \Delta, B, A,
      \Pi}{\Gamma \ENTAIL \Delta, A, B, \Pi} \end{equation}

  The logical rules are those of classical sequent calculus with a
  modified $R\neg$-rule.  The logical rules for negation are used to
  negate atomic formulae only, since the language of $S^i_2E$ allows
  negation of atomic formulae only.

\paragraph{$\neg$-rules:}
\begin{equation}
  \infer[\LEFT\neg]{\neg p(t_1, t_2),
    \Gamma \ENTAIL \Delta}{\Gamma \ENTAIL \Delta, p(t_1, t_2)}
\end{equation}\begin{equation}\label{eq:Rneg}
  \infer[\RIGHT\neg,]{Et_1,Et_2, \Gamma \ENTAIL \Delta, \neg p(t_1,
  t_2)}{p(t_1,
    t_2), \Gamma \ENTAIL \Delta}
\end{equation} where $p$ is $=$ or $\leq$ and $t_1,t_2$ are terms.
Unlike the usual textbook definition of the $\RIGHT\neg$-rule, in the
antecedent we introduce the formulae $Et_1$ and $Et_2$ that express
the existence of values of the terms $t_1$ and $t_2$,
respectively.  This is because we
interpret $\neg p(t_1, t_2)$ to mean that the values $v_1, v_2$ of
$t_1$, $t_2$ exist and satisfy $\neg p(v_1, v_2)$.  Therefore, to
infer $\neg p(t_1,t_2)$ we presuppose $Et_1$ and $Et_2$.

\paragraph{$\AND$-rules:}\begin{equation} \infer[\LEFT\AND_1]{A \AND
B, \Gamma
    \ENTAIL \Delta}{A, \Gamma \ENTAIL \Delta} \end{equation}
  \begin{equation} \infer[\LEFT\AND_2]{B \AND A, \Gamma \ENTAIL
  \Delta}{A, \Gamma
    \ENTAIL \Delta} \end{equation}
  \begin{equation} \infer[\RIGHT\AND]{\Gamma \ENTAIL \Delta, A \AND
    B}{\Gamma \ENTAIL \Delta, A \quad \Gamma \ENTAIL \Delta,
    B} \end{equation}

\paragraph{$\OR$-rules:}\begin{equation} \infer[\LEFT\OR]{A \OR B,
\Gamma \ENTAIL
    \Delta}{A, \Gamma \ENTAIL \Delta \quad B, \Gamma \ENTAIL
    \Delta} \end{equation}
  \begin{equation} \infer[\RIGHT\OR_1]{\Gamma \ENTAIL \Delta, A \OR
    B}{\Gamma \ENTAIL \Delta, A} \end{equation}
  \begin{equation} \infer[\RIGHT\OR_2]{\Gamma \ENTAIL \Delta, B \OR
    A}{\Gamma \ENTAIL \Delta, A} \end{equation}

The logical rules for universal quantification are of two types:
\emph{bounded} (indicated with a ``b'' following the L or R) and
\emph{unbounded}, corresponding to bounded and unbounded universal
quantification, respectively, in the language of $S^i_2E$.  The same
is true of the logical rules for existential quantification. In both
types, bounded and unbounded, the quantification (universal or
existential) is done over objects that actually exist.

\paragraph{Bounded $\forall$-rules:}
\begin{equation}
    \infer[\LEFTb\forall,]{t \leq s, \forall x \leq s.
      A(x), \Gamma \ENTAIL \Delta}{A(t), \Gamma \ENTAIL \Delta}
  \end{equation}where the variable $x$ does not occur in the term $s$.
\begin{equation}
    \infer[\RIGHTb\forall,]{Et, \Gamma \ENTAIL \Delta, \forall
      x \leq t. A(x)}{a \leq t, \Gamma \ENTAIL \Delta, A(a)}
  \end{equation} where neither the variable $a$ nor the
  variable $x$ occurs in the term $t$, and $a$ does not occur free in
  $\Gamma\ENTAIL\Delta$. 

\paragraph{Unbounded $\forall$-rules:}
\begin{equation}
    \infer[\LEFT\forall]{Et, \forall x.
      A(x), \Gamma \ENTAIL \Delta}{A(t), \Gamma \ENTAIL \Delta}
  \end{equation}
  \begin{equation}
    \infer[\RIGHT\forall,]{\Gamma \ENTAIL \Delta, \forall x.
      A(x)}{Ea, \Gamma \ENTAIL \Delta, A(a)}
  \end{equation} where the variable $a$ does not occur free in
  $\Gamma\ENTAIL\Delta$.

\paragraph{Bounded $\exists$-rules:}\begin{equation}
    \infer[\LEFTb\exists,]{\exists x \leq t. A(x),
      \Gamma \ENTAIL \Delta}{a \leq t, A(a), \Gamma \ENTAIL \Delta}
  \end{equation} where neither the variable $a$ nor the
  variable $x$ occurs in neither the term  $t$, and $a$ does not occur free in
  $\Gamma\ENTAIL\Delta$.
 \begin{equation}
    \infer[\RIGHTb\exists,]{t \leq s, \Gamma \ENTAIL \Delta,
      \exists x \leq s. A(x)}{\Gamma \ENTAIL \Delta, A(t)}
  \end{equation}where the variable $x$
  does not occur in the term $s$.

\paragraph{Unbounded $\exists$-rules:}
\begin{equation} \infer[\LEFT\exists,]{\exists x. A(x),
      \Gamma \ENTAIL \Delta}{Ea, A(a), \Gamma \ENTAIL \Delta}
  \end{equation} where the variable $a$ does not occur free in
  $\Gamma\ENTAIL\Delta$.

\begin{equation}
  \infer[\RIGHT\exists]{Et, \Gamma \ENTAIL \Delta,
    \exists x. A(x)}{\Gamma \ENTAIL \Delta, A(t)}
\end{equation}

The \emph{Cut rule} is used to derive a sequent by employing an
intermediate ``lemma.''

\paragraph{Cut rule:}
 \begin{equation}\label{eq:Cut}
   \infer[\textup{Cut}]{\Gamma, \Pi \ENTAIL \Delta,
     \Lambda}{\Gamma \ENTAIL \Delta, A \quad A, \Pi \ENTAIL \Lambda}
 \end{equation}

 The \emph{$\PIND$ rule} is used to infer statements $A(t)$ (for
 formulae $A(a)$ and terms $t$) from the assumption that $t$ converges
 and $A$ satisfies the induction hypothesis on the value of $t$.

\paragraph{$\PIND$ rule:}\begin{equation}
  \label{eq:PIND} \infer[\PIND,]{Et, \Gamma \ENTAIL \Delta, A(t)}{
    \Gamma \ENTAIL \Delta, A(0) \quad A(a), \Gamma \ENTAIL \Delta,
    A(s_0a) \quad A(a), \Gamma \ENTAIL \Delta, A(s_1a)}\end{equation}
where the variable $a$ does not occur free in
$\Gamma\ENTAIL\Delta$.  We call this the $\PIND$ rule on $t$. If we
restrict the $\PIND$ rule to formulae $A(a)$ from some class $\Phi$,
we call the resulting rule the $\Phi\text{-}\PIND$ rule.

Informally, the $\PIND$ rule expresses induction on the notations of
the natural numbers represented in binary.  We call $A(a)$ the
\emph{induction hypothesis}, the proof of $\Gamma \ENTAIL \Delta,
A(0)$ the \emph{base case}, and the proofs of $A(a), \Gamma \ENTAIL
\Delta, A(s_0a)$ and $A(a), \Gamma \ENTAIL \Delta, A(s_1a)$ the
\emph{induction steps}.

From the $\PIND$ rule on $t$, we have the following derived rule.
\begin{equation}
  \label{eq:PIND-E} \infer{Et, \Gamma \ENTAIL \Delta, A(t)}{
    \Gamma \ENTAIL \Delta, A(0) \quad Ea, A(a), \Gamma \ENTAIL \Delta,
    A(s_0a) \quad Ea, A(a), \Gamma \ENTAIL \Delta, A(s_1a)
  }
\end{equation} We call this the $\PIND\text{-}E$ rule on
$t$.  To derive this rule, we first derive $A(a) \AND Ea, \Gamma
\ENTAIL \Delta, A(s_0a) \AND Es_0a$ from $Ea, A(a), \Gamma \ENTAIL
\Delta, A(s_0a)$.  This is done via the axiom $Ea \ENTAIL Es_0a$ and
propositional reasoning.  Similarly, we derive $A(a) \AND Ea, \Gamma
\ENTAIL \Delta, A(s_1a) \AND Es_1a$ from $Ea, A(a), \Gamma \ENTAIL
\Delta, A(s_1a)$.  Then, applying $\PIND$ to $A(a) \AND Ea$, we derive
$Et, \Gamma \ENTAIL \Delta, A(t) \AND Et$.  Finally, we derive $Et,
\Gamma \ENTAIL \Delta, A(t)$ by purely propositional reasoning.

This completes the definition of the rules of inference of $S^i_2E$.

\begin{lem}[Substitution Lemma]\label{lem:subst}
  If $S^i_2E\vdash \Gamma(a) \ENTAIL \Delta(a)$, where the variable
  $a$ occurs free in $\Gamma\ENTAIL\Delta$, then
  $S^i_2E\vdash\Gamma(t)\ENTAIL\Delta(t)$ for every term $t$.
\end{lem}

\proof
  Induction on the structure of the $S^i_2E$ proof of $\Gamma(a)
  \ENTAIL \Delta(a)$.
\qed

\subsection{Comparison to Free Logic and Scott's system}
\label{subsec:freelogic}

In this subsection, we briefly discuss the similarities between our
system and two others: free logic (see \cite{FreeLogic} for a
comprehensive summary) and Scott's $E$-logic \cite{Scott1979}.

Our system is quite similar to free logic with negative semantics,
known as NFL, and Scott's system, in that in all three systems an
equational formula $t=u$ implies existence of the objects denoted by
$t$ and $u$.  Moreover, in our system as well as in NFL, all the
functions are \emph{strict} (in Scott's terminology), that is,
$Ef(a_1, \dots, a_n)$ implies $Ea_j$ ($j = 1, \dots, n$).  The rules
for universal and existential quantification are also similar.

On the other hand, our $R\neg$ rule is more restrictive than that of
either NFL or Scott's system.  Also, we do not have \emph{full
  preservation of equality} for functions, $a_1 = b_1, \dots, a_n =
b_n \ENTAIL f(a_1, \dots, a_n) = f(b_1, \dots, b_n)$, as axioms
(although that can be proved as a theorem).  Furthermore, we do not
have \emph{totality} of functions, $Ea_1, \dots, Ea_n \ENTAIL Ef(a_1,
\dots, a_n)$, as axioms (again, that can be proved as a theorem).  The
only functions for which we have totality as axioms are the binary
successor functions $s_0$ and $s_1$; however, we can prove the
totality of the other functions from the binary successor functions
via the $\Sigma^b_0\text{-}\PIND$ rule.

As we will see in the next section, the restrictions we have imposed
on our system are crucial to the provability in $S_2$ of the
consistency of strictly $i$-normal proofs in $S^{-1}_2E$.

\section{$S^{i+2}_2$ proof of consistency of strictly $i$-normal
  proofs}\label{sec:consistency}

In this section, we define \emph{$i$-normal formula} and
\emph{strictly $i$-normal proof}, and in $S^{i+2}_2$ we prove the
consistency of strictly $i$-normal proofs in $S^{-1}_2E$.  The
consistency proof is based on the facts that we can produce a
$\Sigma^b_i$ formula that constitutes a truth definition for
$i$-normal formulae and we can apply the $\Sigma^b_{i+2}\text{-}\PIND$
rule to prove the soundness of strictly $i$-normal proofs in
$S^{-1}_2E$.  The idea is that to use a term $t$ in an $S^{-1}_2E$
proof, we first need to prove that $Et$ holds.  To do this, we show
that for a given assignment $\rho$ of values to the variables in $t$,
the value of $t$ is bounded by the size of the proof of $Et$ plus the
size of $\rho$.  Therefore, we can define a valuation function for
terms and a truth definition for the formulae in the proof.  Once we
obtain the truth definition, consistency is easy to prove.

In Subsection \ref{subsec:goedel}, we describe how we assign G\"odel
numbers to formulae and proofs.  In Subsection \ref{subsec:value}, we
introduce a ``bounded'' valuation for terms and prove its basic
properties.  In Subsection \ref{subsec:truth}, we introduce a
``bounded'' truth definition for quantifier-free formulae and then
extend it to $i$-normal formulae.  In Subsection \ref{subsec:sound},
we introduce strictly $i$-normal proofs in $S^{-1}_2E$, and we prove
the soundness of such proofs with respect to the truth definition
given in Subsection \ref{subsec:truth}.

In what follows, $i$ is assumed to be a fixed integer from the set
$\{-1, 0,1,2,\dots\}$.

\subsection{G\"odel numbers for formulae and
proofs }\label{subsec:goedel} Before proceeding to
the main topics treated in this section, we need to explain how we
assign G\"odel numbers to the formulae of $S^i_2E$ and to $S^i_2E$
proofs.

Since our proofs are structured as proof trees (nested sequences of
formulae and rules of inference), we first assign G\"odel numbers to
the symbols which appear in the formulae and rules of inference of
$S^i_2E$, using a different natural number for each symbol.  Then we
assign G\"odel numbers to finite sequences of natural numbers, and we
apply this G\"odel numbering to the formulae of $S^i_2E$ and the nodes
of proof trees.

In assigning G\"odel numbers to sequences of natural numbers, we
follow Buss' method \cite{Buss:Book}.  The G\"odel
number of the sequence $\langle u_1, \ldots, u_n \rangle$ is
determined by the following procedure.
\begin{enumerate}[(1)]
\item For every $j$ ($1\le j\le n$), let $u_j'$ be the natural number
  whose binary representation is obtained by inserting 0 between
  consecutive bits in the binary representation of $u_j$, and then
  appending 0 to the resulting bit string. (By convention, $0'=0$.)
  For example, if the binary representation of $u_j$ is 11, then that
  of $u_j'$ is 1010; and if the binary representation of $u_j$ is 101,
  then that of $u'$ is 100010.
\item Define the G\"odel number of $\langle u_1, \ldots, u_n \rangle$
  to be $u_1'\oplus 3\oplus u_2'\oplus 3\oplus\cdots \oplus 3\oplus
  u_{n-1}' 3\oplus u_n'$.  Since the binary representation of the
  number 3 is 11, the G\"odel number of the sequence $\langle
  u_1,\ldots,u_n \rangle$ is the natural number whose binary
  representation is obtained by inserting 11 between the binary
  representations of $u_j'$ and $u_{j+1}'$ ($j=1,\ldots,n-1$).
\item The G\"odel number of $\langle u_1 \rangle$ is defined to be
  $u_1'$, and the G\"odel number of the empty sequence $\langle
  \rangle$ defined to be 0.
\end{enumerate} The notation $\lceil\cdot\rceil$ is used for the
G\"odel number of any symbol or sequence. For example, the G\"odel
number of the formula $A(a)$ is denoted by $\lceil A(a)\rceil$.

\subsection{Valuation of terms inside $S^i_2$}
\label{subsec:value}

In this subsection, we define a valuation function for the terms of
$S^i_2E$.  Our strategy in defining this valuation is to attach values
to the nodes of a tree which is made up of all the subterms of the
given term, and then to define the value of the term as the value
attached to the root of the tree (the node that represents the entire
term).  We can view construction of this tree as a process of
computation we undertake to obtain the value of its root.

\begin{defi}\label{defn:env}
  Let $\rho$ be a finite sequence of pairs $(\lceil x\rceil,n)$
  where $x$ is a variable of $S^i_2E$ and $n\in\mathbb{N}$.
\begin{enumerate}[(1)]
\item $\rho$ is an \emph{environment} if, for every variable $x$ of
  $S^i_2E$, there is at most one $n\in\mathbb{N}$ with $(\lceil
  x\rceil,n)\in\rho$.
\item If $\rho$ is an environment, we write $\rho(x)=n$ to denote that
  $(\lceil x\rceil,n)\in\rho$, and we sometimes write $\rho(\vec{a})$
  to denote the sequence $\rho(a_1), \dots, \rho(a_k)$ for a finite
  sequence $\vec{a}\equiv a_1,\dots,a_k$ of variables of
  $S^i_2E$.
\item Let $\sigma$ be a term, a formula, or a sequent. $\rho$ is an
  \emph{environment for $\sigma$} if $\rho$ is an environment and, for
  every variable $x$ that occurs free in $\sigma$, there is a pair
  $(\lceil x\rceil,n)$ in $\rho$ (for some $n\in\mathbb{N}$).
\item If $\rho$ is an environment, and if $x$ is a variable and
  $n\in\mathbb{N}$, then $\rho[x \mapsto n]$ denotes environment
  obtained from $\rho$ by replacing the pair $(\lceil
  x\rceil,\rho(x))$ with $(\lceil x \rceil, n)$ if there
  is some $m\in\mathbb{N}$ with $(\lceil x\rceil, m)$ in $\rho$, and
  by adding the pair $(\lceil x\rceil,n)$ to $\rho$ otherwise.
\item $\Env$ denotes the ternary relation that holds of precisely the
  triples $(\rho',\lceil\sigma\rceil, u$) where $\sigma$ is a term, a
  formula, or a sequent; $\rho'$ is an environment for $\sigma$;
  $u\in\mathbb{N}$; for every variable $x$ of $S^i_2E$, there is a
  pair $(\lceil x\rceil, n)$ in $\rho'$ (for some $n\in\mathbb{N}$) if
  and only if $x$ occurs free in $\sigma$; and $\rho'(x)\leq u$ for
  every variable $x$ that occurs free in $\sigma$.  From this point
  on, we identify environments with their G\"odel numbers; therefore,
  we regard $\Env$ as a ternary relation on $\mathbb{N}$.
\item Let $u\in\mathbb{N}$, and let $\sigma$ be a term, a formula, or
  a sequent.  $\BdEnv(\lceil \sigma \rceil, u)$ denotes the greatest
$m\in\mathbb{N}$ which is (the G\"odel number of) an environment
$\rho'$ such that $\Env(\rho',\lceil \sigma \rceil, u)$ holds.
\end{enumerate} 
Although an environment is a sequence---not a set---of pairs, from
this point on we make no reference to the order of the pairs in an
environment.
\end{defi}

For a term, formula, or sequent $\sigma$ and natural numbers $\rho, u$, the relation
$\Env(\rho,\lceil\sigma\rceil,u)$ and the function
$\BdEnv(\lceil\sigma\rceil,u)$ are $\Sigma^b_1$-definable. Also, if
$t$ is a term and $\rho$ is an environment for $t$, we can extend
$\rho$ to $t$ by recursion on the construction of $t$, as
follows: \begin{align} \rho(0) &:= 0 \label{eq:rho_t_0}\\ \rho(f(t_1,
  \dotsc, t_k)) &:= f(\rho(t_1), \dotsc, \rho(t_k))
\end{align} For a fixed
term $t$, the function that, given an arbitrary environment $\rho'$
for $t$, maps $\rho'$ to the value of $\rho'(t)$ obtained by this
method is definable in $S^1_2$.  For a fixed environment $\rho$,
however, the function that, given an arbitrary term $t'$ such that
$\rho$ is an environment for $t'$, maps $t'$ to the value of
$\rho(t')$ obtained by this method is not definable in $S_2$, since
there are sequences of such terms
$t'$ for which the values of $\rho(t')$ increase exponentially in the
length of $t'$. 

\begin{defi}
  Let $t$ be a term of $S^i_2E$, let $\rho$ be an environment for $t$,
  and let $u\in\mathbb{N}$. A \emph{$\rho$-valuation tree for $t$
    which is bounded by $u$} is a tree $w$ that satisfies the
  following conditions.
\begin{enumerate}[(1)]

\item Every node of $w$ is of the form $\langle \lceil t_j\rceil,c
  \rangle$ where $t_j$ is a subterm of $t$, $\;c\in\mathbb{N}$, and
  $c\leq u$.

  \item Every leaf of $w$ is either $\langle \lceil 0 \rceil, 0
    \rangle$ or $\langle \lceil a \rceil, \rho(a)
    \rangle$ for some variable $a$ in the domain of $\rho$.

  \item The root of $w$ is $\langle \lceil t \rceil, c
    \rangle$ for some $c\leq u$.

  \item If $\langle \lceil f(t_1, \dots, t_n) \rceil, c \rangle$ is a
    node of $w$, then the children of this node are the nodes $\langle
    \lceil t_1 \rceil, d_1 \rangle, \dots, \langle \lceil t_n \rceil,
    d_n \rangle$ which satisfy the condition $c = f(d_1, \dots, d_n)$.
  \end{enumerate}
\end{defi}

If the root of a $\rho$-valuation tree $w$ for $t$ is $\langle
\lceil t \rceil, c \rangle$, we say \emph{the value of $w$ is $c$}.

The statement that $t$ converges to the value $c$ (and $c\le u$) is
defined by the $\Sigma^b_1$ formula which expresses that the following
relation (which we denote by $v(\lceil t \rceil, \rho) \downarrow_u
c$) holds: ``$ \exists w \leq s(\lceil t \rceil, u)$ such that $w$ is
(the G\"odel number of) a $\rho$-valuation tree for $t$ which is
bounded by $u$ and has root $\langle \lceil t \rceil, c \rangle$,''
where $s(\lceil t\rceil,u)$ is a term whose G\"odel number bounds (the
G\"odel numbers of) all $\rho$-valuation trees for $t$ which are
bounded by $u$.

Here are some simple facts about the relation $v(\lceil t \rceil,
\rho) \downarrow_u c$.

\begin{lem}\label{lem:v} For a term $t$, an environment $\rho$ for
$t$, and $u\in\mathbb{N}$, the following statements are provable in
$S^1_2$.
 \begin{enumerate}[\em(1)]
 \item If $v(\lceil t \rceil, \rho) \downarrow_u c$ and $v(\lceil t
   \rceil, \rho) \downarrow_u c'$, then $c = c'$.
 \item If $f$ is an $n$-ary function such that $t\equiv ft_1\cdots
   t_n$ for terms $t_1,\dots,t_n$ and we have $v(\lceil f(t_1, \dots, t_n)
   \rceil, \rho) \downarrow_u c$, then there exist $d_1, \dots, d_n$
   such that $f(d_1, \dots,
   d_n) = c$ and $v(\lceil t_1 \rceil, \rho) \downarrow_u d_1,\ldots,
   v(\lceil t_n \rceil, \rho) \downarrow_u d_n$.
 \item $v(\lceil 0 \rceil, \rho) \downarrow_0 0$
\item For every variable $a$ that occurs in $t$, $v(\lceil a
   \rceil,\rho)\downarrow_{\rho(a)}\rho(a)$.
 \item If $v(\lceil t \rceil, \rho) \downarrow_u c$ and $\rho' \subset
   \rho$ is an environment for $t$, $v(\lceil t \rceil, \rho')
   \downarrow_u c'$ and $c = c'$.
 \end{enumerate}
\end{lem}

Lemmata \ref{lem:v-sub}, \ref{lem:v-upward} and
\ref{lem:v-standard} are key properties of $v$ that are used in our
proof of the consistency of strictly $i$-normal proofs
(Subsection~\ref{subsec:sound}).  Lemma \ref{lem:v-sub} states that
substitution is provably equivalent to assignment.  Lemma
\ref{lem:v-upward} states that the relation $v(\lceil t \rceil, \rho)
\downarrow_u c$ is closed upward with respect to $u$.  Lemma
\ref{lem:v-standard} states that if $u$ is sufficiently large, the
relation $v(\lceil t \rceil, \rho) \downarrow_u c$ coincides with the
valuation of $t$.

\begin{lem}\label{lem:v-sub}
  For terms $t(\vec{a},a)$ and $t'(\vec{a})$ such that the
variable $a$ occurs in $t$ and is not in $\vec{a}$, the following is
provable in $S^1_2$.
\[ v(\lceil t(\vec{a}, t'(\vec{a})) \rceil, \rho) \downarrow_u c
\leftrightarrow \exists c' \leq u, v(\lceil t'(\vec{a}) \rceil, \rho)
\downarrow_u c' \AND v(\lceil t(\vec{a}, a) \rceil, \rho[a \mapsto
c']) \downarrow_u c\]
\end{lem}

\proof

By induction on the contraction of $t$.

\qed

\begin{lem}\label{lem:v-upward} The following statement is
  provable in $S^1_2$: ``If $v(\lceil t \rceil, \rho) \downarrow_u c$
  and $u < u'$, then $v(\lceil t \rceil, \rho) \downarrow_{u'} c$.''
\end{lem}

\proof

The lemma holds since if $\rho$-valuation tree $w$ is bounded by $u$,
then it is bounded by $u'$.

\qed

\begin{lem}\label{lem:v-standard}
  Let $t$ be a term of $S^i_2E$, let $t_1, \dots, t_m$ be an
  enumeration of all the subterms of $t$.  Then the following
  statement is provable in $S^1_2$: ``For any environment $\rho$ for
  $t$ and $u\in\mathbb{N}$, if $\rho(t_j)\leq u$ for every $j$, then
  $v(\lceil t \rceil, \rho) \downarrow_u \rho(t)$ holds; and if
  $\rho(t_j)>u$ for some $j$, then $v(\lceil t \rceil,
  \rho)\downarrow_u c$ does not hold for any natural number $c\leq
  u$.''.
\end{lem}

\proof
  If $\rho(t_j)\leq u$ for every $j$, we can construct a
  $\rho$-valuation tree for $t$ which is bounded by $u$, by induction
  (outside of $S^1_2$) on $t$.  If $\rho(t_j)>u$ for some $j$,
  then there is no $\rho$-valuation tree for $t$ which is
  bounded by $u$ and this fact can be proved in $S^1_2$.
\qed

\subsection{Truth definition inside $S^i_2$}
\label{subsec:truth}

In this subsection, we give a ``bounded'' truth definition for
quantifier-free formulae, and then we extend the definition to the
formulae $i$-normal formulae \cite{Takeuti:Truth}, \cite{Takeuti1996},
\cite{Takeuti2000}.

First, we present a truth definition for quantifier-free formulae.
Since logical symbols can be arbitrarily nested, we follow the same
strategy that was used in our definition of valuation for terms.  We
attach a truth value to each node of a subformula tree, and we define
the value attached to the root (the node that represents the entire
formula) as the truth value of the formula.

\begin{defi}\label{defn:T_{-1}}
  Let $A$ be a quantifier-free formula of $S^{-1}_2E$, let $\rho$ be
  an environment for $A$, and let $u\in\mathbb{N}$. A
  \emph{$\rho$-truth tree for $A$ which is
    bounded by $u$} is a tree $w$ that satisfies the following
  conditions.

Every leaf of $w$ has one of the following five forms (where in each
form the possible values of $\epsilon$ are $0$ and $1$): $\langle
\lceil t_1 \leq t_2 \rceil, \epsilon \rangle$, $\langle \lceil t_1
\not\leq t_2 \rceil, \epsilon \rangle$, $\langle \lceil t_1 = t_2
\rceil, \epsilon \rangle$, $\langle \lceil t_1 \not= t_2 \rceil,
\epsilon \rangle$, $\langle \lceil Et \rceil, \epsilon \rangle$.

For a leaf of the form $\langle \lceil t_1 \leq t_2 \rceil, \epsilon
\rangle$, $\epsilon = 1$ if $\exists c_1, c_2 \leq u$, $v(\lceil t_1
\rceil, \rho) \downarrow_u c_1$, $v(\lceil t_2 \rceil, \rho)
\downarrow_u c_2$, and $c_1 \leq c_2$; otherwise, $\epsilon = 0$.

For a leaf of the form $\langle \lceil t_1 \not\leq t_2 \rceil,
\epsilon \rangle$, $\epsilon = 1$ if $\exists c_1, c_2 \leq u$,
$v(\lceil t_1 \rceil, \rho) \downarrow_u c_1$, $v(\lceil t_2 \rceil,
\rho) \downarrow_u c_2$, and $c_1 \not\leq c_2$; otherwise, $\epsilon
= 0$.

The conditions that must be satisfied by a leaf of the form $\langle
\lceil t_1 = t_2 \rceil, \epsilon \rangle$ or $\langle \lceil t_1
\not= t_2 \rceil, \epsilon \rangle$ are the obvious analogues of those
for $\langle \lceil t_1 \leq t_2 \rceil, \epsilon \rangle$ and
$\langle \lceil t_1 \not\leq t_2 \rceil, \epsilon \rangle$,
respectively.

For a leaf of the form $\langle \lceil Et \rceil, \epsilon \rangle$,
$\epsilon = 1$ if $\exists c \leq u$, $v(\lceil t \rceil, \rho)
\downarrow_u c$; otherwise, $\epsilon = 0$.

Every intermediate node $r$ of $w$ is of the form $\langle \lceil A_1
\AND A_2 \rceil, \epsilon \rangle$ or $\langle \lceil A_1 \OR A_2
\rceil, \epsilon \rangle$, where the children of $r$ are the nodes
$\langle \lceil A_1 \rceil, \epsilon_1 \rangle$ and $\langle \lceil
A_2 \rceil, \epsilon_2 \rangle$.

For a node of the form $\langle \lceil A_1 \AND A_2 \rceil, \epsilon
\rangle$, $\epsilon = 1$ if $\epsilon_1 = 1$ and $\epsilon_2 = 1$;
otherwise, $\epsilon = 0$.

For a node of the form $\langle \lceil A_1 \OR A_2 \rceil, \epsilon
\rangle$, $\epsilon = 1$ if $\epsilon_1 = 1$ or $\epsilon_2 = 1$;
otherwise, $\epsilon = 0$.

The root of $w$ is $\langle \lceil A \rceil, \epsilon \rangle$ for
some $\epsilon\in\{0,1\}$.

The truth of a quantifier-free formula $A$ is
defined by the $\Sigma^b_1$ formula $T_{-1}(u, \lceil A \rceil, \rho)$
which expresses that ``$\exists w \leq s(\lceil A \rceil, u)$ such
that $w$ is (the G\"odel number of) a $\rho$-truth tree for $A$ which
is bounded by $u$ and has root $\langle
\lceil A \rceil, 1 \rangle$,'' where $s(\lceil A \rceil, u)$ is
a term which bounds (the G\"odel numbers of) all $\rho$-truth trees for $A$
which are bounded by $u$.
\end{defi}

We can prove several basic properties of the truth definition
$T_{-1}$.
\begin{lem}\label{lem:T_{-1}}
  The following statements are provable in $S^1_2$.
\begin{enumerate}[\em(1)]
 \item $T_{-1}(u, \lceil t_1 \leq t_2 \rceil, \rho) \leftrightarrow
 \exists
       c_1, c_2 \leq u, v(\lceil t_1 \rceil, \rho) \downarrow_u c_1
       \AND v(\lceil t_2 \rceil, \rho) \downarrow_u c_2 \AND c_1 \leq
       c_2$
 \item $T_{-1}(u, \lceil t_1 \not\leq t_2 \rceil, \rho)
 \leftrightarrow \exists
       c_1, c_2 \leq u, v(\lceil t_1 \rceil, \rho) \downarrow_u c_1
       \AND v(\lceil t_2 \rceil, \rho) \downarrow_u c_2 \AND c_1
       \not\leq c_2$
 \item $T_{-1}(u, \lceil t_1 = t_2 \rceil, \rho) \leftrightarrow
 \exists
       c_1, c_2 \leq u, v(\lceil t_1 \rceil, \rho) \downarrow_u c_1
       \AND v(\lceil t_2 \rceil, \rho) \downarrow_u c_2 \AND c_1 =
       c_2$
 \item $T_{-1}(u, \lceil t_1 \not= t_2 \rceil, \rho) \leftrightarrow
 \exists
       c_1, c_2 \leq u, v(\lceil t_1 \rceil, \rho) \downarrow_u c_1
       \AND v(\lceil t_2 \rceil, \rho) \downarrow_u c_2 \AND c_1 \not=
       c_2$
 \item $T_{-1}(u, \lceil Et \rceil, \rho) \leftrightarrow \exists
       c \leq u, v(\lceil t \rceil, \rho) \downarrow_u c$
 \item $T_{-1}(u, \lceil A_1 \AND A_2 \rceil, \rho)
       \leftrightarrow T_{-1}(u, \lceil A_1 \rceil, \rho) \AND
       T_{-1}(u, \lceil A_2 \rceil, \rho)$
 \item $T_{-1}(u, \lceil A_1 \OR A_2 \rceil, \rho)
       \leftrightarrow T_{-1}(u, \lceil A_1 \rceil, \rho) \OR
       T_{-1}(u, \lceil A_2 \rceil, \rho)$
 \item $T_{-1}(u, \lceil A(\vec{a}, t(\vec{a})) \rceil, \rho)
       \leftrightarrow \exists c \leq u, v(\lceil t(\vec{a}) \rceil,
       \rho) \downarrow_u c \AND T_{-1}(u, \lceil A(\vec{a}, a)
       \rceil, \rho[a \mapsto c])$
     \item If $\rho' \subset \rho$ is an environment for $A$,
       $T_{-1}(u, \lceil A \rceil, \rho') \leftrightarrow T_{-1}(u,
       \lceil A \rceil, \rho)$
\end{enumerate}
\end{lem}

Lemmata \ref{lem:T0-sub}, \ref{lem:T0-upward}, \ref{lem:T0-refl}, and
\ref{lem:EM} are key properties of $T_{-1}$ that are used in our
consistency proof.  Lemma \ref{lem:T0-sub} states that substitution is
provably equivalent to assignment.  Lemma \ref{lem:T0-upward} states
that $T_{-1}$ is closed upward with respect to $u$.  Lemma
\ref{lem:T0-refl} states a reflection principle for $T_{-1}$. Lemma
\ref{lem:EM} states the law of the excluded middle.

\begin{lem}\label{lem:T0-sub}
  Let $A$ be a quantifier-free formula in which the variable $a$
  occurs, and let $t$ be a term. Then the following statement is
  provable in $S^1_2$:
  ``$T_{-1}(u, \lceil A(t) \rceil, \rho)$ if and only if there exists
  $c\leq u$ such that $v(\lceil t \rceil,\rho) \downarrow_u c$ and
  $T_{-1}(u, \lceil A(a) \rceil, \rho[a \mapsto c])$.''
\end{lem}

\begin{lem}\label{lem:T0-upward}
  For a quantifier-free formula $A$, it is provable in $S^1_2$ that
  $T_{-1}(u, \lceil A \rceil, \rho), u < u' \rightarrow T_{-1}(u',
  \lceil A \rceil, \rho)$.
\end{lem}

\proof
  Induction on the construction of $A$, using Lemma \ref{lem:v-upward}
  for basic formulae $A$, and clauses 6 and 7 of Lemma
  \ref{lem:T_{-1}} for other quantifier-free formulae.
\qed

We can prove that $T_{-1}$ is a truth definition by showing a kind of
reflection principle for $T_{-1}$.
\begin{lem}\label{lem:T0-refl} Let $A(\vec{a})$ be a
  quantifier-free formula, let $A'(\vec{a})$ be the interpretation of
  $A(\vec{a})$ which is given in Definition \ref{defn:interpretation}
  and let $t_1, \dots, t_m$ be an enumeration of all the subterms of
  the terms that occur in $A(\vec{a})$.  Then the following statement
  is provable in $S^1_2$: ``If $\rho$ be an environment for $A$,
  $u\in\mathbb{N}$ and $\rho(t_j) \leq u$ for every $j$, then
  $T_{-1}(u, \lceil A(\vec{a}) \rceil, \rho) \leftrightarrow
  A'(\rho(\vec{a}))$; and if $\rho(t_j)>u$ for some
  $j$, then $\neg T_{-1}(u, \lceil A(\vec{a}) \rceil, \rho)$.''
\end{lem}

\proof
  Induction on the construction of $A$, using Lemma
  \ref{lem:v-standard} for basic formulae $A$.
\qed

\begin{lem}\label{lem:EM} If $p$ is $=$ or $\leq$, it is provable in
  $S^1_2$ that if $T_{-1}(u, \lceil Et_i \rceil, \rho)$ for  $i \in \{1, 2\}$ hold, then either $T_{-1}(u, \lceil
  p(t_1, t_2) \rceil, \rho)$ or $T_{-1}(u, \lceil \neg p(t_1, t_2)
  \rceil, \rho)$ holds.
\end{lem}

Next, we would like to present a truth definition for $\Sigma^b_i$
formulae.  However, since it is technically difficult to do this for
general $\Sigma^b_i$ formulae, we restrict our definition to
\emph{$i$-normal} formulae.  Since $i\in\{-1,0,1,2,\ldots\}$, we have
$-1$-normal formulae, $0$-normal formulae, $1$-normal formulae,
$2$-normal formulae, and so on.
\begin{defi}\label{defn:normal}
 Let $i\ge -1$, and let $A(\vec{a})$ be a formula.

If $i=-1$, $A(\vec{a})$ is \emph{pure $-1$-normal} if
$A(\vec{a})$ is quantifier free.

If $i\ge 0$, $A(\vec{a})$ is \emph{pure
$i$-normal} if it is of the
  form

\begin{multline*} \exists x_1 \leq t_1(\vec{a}) \forall x_2 \leq
t_2(\vec{a},
 x_1)\cdots \\ \qquad Q_i x_i \leq
t_i(\vec{a},x_1,\dots,x_{i-1})Q_{i+1} x_{i+1} \leq
|t_{i+1}(\vec{a},x_1,\dots,x_i)|.A_0(\vec{a},x_1,\dots,x_{i+1}),
\end{multline*} 
where $Q_i$ is $\forall$ if $i$ is even, and $\exists$ if $i$ is odd
(and vice versa for $Q_{i+1}$), and $A_0(\vec{a},x_1,\dots,x_{i+1})$
is quantifier free and does not contain the predicate $E$.  Note that
if $i\ge 0$, then a pure $i$-normal formula has at least one
quantifier, so quantifier-free formulae are not pure $i$-normal.  Also
note that unlike the variables $x_1,\ldots,x_i$ (which are bounded by
the values of the terms $t_1,\ldots,t_i$, respectively), the value of
the variable $x_{i+1}$ is bounded by
$|t_{i+1}(\vec{a},x_1,\dots,x_i)|$, that is, by the number of bits in
the binary representation of the value of the term
$t_{i+1}(\vec{a},x_1,\dots,x_i)$.

If $i=-1$, $A(\vec{a})$ is \emph{$i$-normal} if it is quantifier
free. 

If $i\ge 0$, $A(\vec{a})$ is \emph{$i$-normal} if it is a subformula
of a pure $i$-normal formula or is $Et$ for some term $t$.  In other
words, $A(\vec{a})$ is either an $E$-form, a
quantifier-free formula that does not contain $E$, or a formula of the
form
\begin{multline}
  Q_jx_j \leq t_j(\vec{a},x_1,\dots,x_{j-1}) \cdots Q_ix_i \leq
  t_i(\vec{a}, x_1, \dots, x_{i-1}) \\
 Q_{i+1} x_{i+1}\leq
 |t_{i+1}(\vec{a},x_1,\dots,x_i)|.A_0(\vec{a},x_1,\dots,x_{i+1}),
\end{multline} where $A_0(\vec{a},x_1,\dots,x_{i+1})$ is quantifier
free and does not contain $E$; $1\leq j\leq i + 1$; and for every $k$
with $j\leq k\leq i+1$, $Q_k$ is either $\forall$ or $\exists$,
according as $k$ is even or odd.  If $j = i + 1$, the above formula is
$Q_{i+1} x_{i+1}\leq
|t_{i+1}(\vec{a},x_1,\dots,x_i)|.A_0(\vec{a},x_1,\dots,x_{i+1})$.
\end{defi}

The following is a truth definition $T_i(u, \lceil B \rceil, \rho)$
for $i$-normal formulae $B$.  First, we define a truth definition
$T_{i, l}$ for $i$-normal forms with $l$ quantifiers.
\begin{defi}\label{defn:T}
  Let $i\ge -1$, let $B$ be an $i$-normal formula with $l$
  quantifiers.  Note that $0\leq l\leq i+1$.  We define $T_{i, l}(u,\lceil
  B\rceil,\rho)$ by recursion on $l$ in the meta-language.

  If $l=0$, then $B$ is quantifier free, so $T_i(u,\lceil
  B\rceil,\rho) \equiv T_{-1}(u,\lceil B\rceil,\rho)$.

  If $l\ge 1$, then
  \[B\equiv Q_jx_j\leq t.A(\vec{a},x_1,\dots,x_{j}),\] where $j=i+2-l$;
  $t\equiv t_j(\vec{a},x_1,\dots,x_{j-1})$ if $j<i+1$, and $t\equiv
  |t_{i+1}(\vec{a},x_1,\dots,x_i)|$ if $j=i+1$; and
  $A(\vec{a},x_1\dots,x_{j})$ is an $i$-normal formula with $l-1$
  quantifiers. Assume that we have defined $T_{i, l-1}(u,\lceil C\rceil,\rho)$
  for all $i$-normal formulae $C$ with $l-1$ quantifiers. We define
  $T_{i, l}(u, \lceil B \rceil, \rho)$ to be the following formula.
  \[\exists c \leq u, v(\lceil t \rceil, \rho)
  \downarrow_u c \AND Q_j d_j \leq c. T_i(u, \lceil A(\vec{a},
  x_1,\dots,x_{j}) \rceil, \rho[x_j \mapsto d_j])\]

  Then, let $\textup{INQ}(\lceil B \rceil, l)$ be a formula which represents
  ``B is an $i$-normal form with $l$ quantifiers''.  we define
  $T_i(u, \lceil B \rceil, \rho)$ as
  \begin{equation}
    \{\textup{INQ}(\lceil B \rceil, 0) \IMPLY T_{i, 0}(u, \lceil B \rceil, \rho)\}
    \OR \ldots \OR \{\textup{INQ}(\lceil B \rceil, i+1) \IMPLY T_{i, i+1}(u,
    \lceil B \rceil, \rho). \}
  \end{equation}

Since we can contract successive $\exists$ quantifiers into a single
$\exists$ quantifier, $T_i(u, \lceil B \rceil, \rho)$ is
$\Sigma^b_{i+1}$.
 \end{defi}

 Lemmata \ref{lem:T-sub} and \ref{lem:T-upward} are key properties of
 $T_i$ that are used in our consistency proof.  Lemma \ref{lem:T-sub}
 states that substitution is provably equivalent to assignment.  Lemma
 \ref{lem:T-upward} states that $T_i$ is closed upward with respect to
 $u$.

\begin{lem}\label{lem:T-sub}
  Let $A$ be an $i$-normal formula in which the variable $a$ occurs,
  and let $t$ be a term. Then the following statement is provable in
  $S^1_2$: ``$T_i(u, \lceil A(t) \rceil, \rho)$ if and only if there
  exists $c\le u$ such that $v(\lceil t \rceil,\rho) \downarrow_u c$
  and $T_i(u, \lceil A(a) \rceil, \rho[a \mapsto c])$.''
\end{lem}

\proof
  Easy consequence of Lemma \ref{lem:T0-sub}.
\qed

Similarly, we can use Lemma \ref{lem:T0-upward} to prove upward
closedness of $T_i$ with respect to $u$.
\begin{lem}\label{lem:T-upward}
  For an $i$-normal formula $A$, it is provable in $S^1_2$ that
  $T_i(u, \lceil A \rceil, \rho) \AND u < u' \rightarrow T_i(u',
  \lceil A \rceil, \rho)$.
\end{lem}

As the clause 9 of Lemma \ref{lem:T_{-1}}, if $\rho$ contains
variables other than free variables which occurs in $A$, we can ignore
such variables.
\begin{lem}\label{lem:T-rho}
  For an $i$-normal formula $A$, environments $\rho$ for $A$ and
  $\rho \subset \rho'$, it is provable in $S^1_2$ that $T_i(u, \lceil
  A \rceil, \rho) \leftrightarrow T_i(u, \lceil A \rceil, \rho')$.
\end{lem}

By definition of $T_i$, we see that we can take the outermost
quantifier of an $i$-normal formula whose G\"odel number is the second
argument of $T_i$ and move it to the outside of $T_i$.
\begin{lem}\label{lem:T}
  Let $\forall x \leq t. A(x)$ and $\exists x \leq t. A(x)$ be
  $i$-normal formulae with at least one quantifier (and whose
  outermost quantifier is of the indicated type).  Then the following
  is provable in $S^1_2$.
  \begin{multline}
    T_i(u, \lceil \forall x \leq t. A(x) \rceil, \rho) \Leftrightarrow
  \\ \exists c \leq u, v(\lceil t \rceil, \rho) \downarrow_u c \AND
  \forall d \leq c. T_i(u, \lceil A(a) \rceil, \rho[a \mapsto d])
  \end{multline}
  \begin{multline} T_i(u, \lceil \exists x \leq t. A(x) \rceil, \rho)
  \Leftrightarrow \\ \exists c \leq u, \ v(\lceil t \rceil, \rho)
  \downarrow_u c \AND \exists d \leq c. T_i(u, \lceil A(a) \rceil,
  \rho[a \mapsto d])
  \end{multline}
\end{lem}

\subsection{Strictly $i$-normal proofs and their consistency}
\label{subsec:sound}

Since we have defined truth for $i$-normal formulae only, we can
define soundness for only those proofs that consist entirely of
$i$-normal formulae.  We call such proofs \emph{strictly $i$-normal},
and we use the term \emph{strictly $i$-normal proof tree for
  $\Gamma\ENTAIL\Delta$} to mean a tree $w$ which represents a proof
of $\Gamma \ENTAIL \Delta$.  In such a tree, every node $r$ has the
form $\langle R, \lceil \Gamma_r \ENTAIL \Delta_r \rceil, w_1, \ldots,
w_{l(R)} \rangle$, where $R$ is a name of inference, $\Gamma_r \ENTAIL
\Delta_r$ is the conclusion of the inference, and $l(R)$ is a number
of premises of $R$.  If $l(R)=0$, the node has the form $\langle R,
\lceil \Gamma_r \ENTAIL \Delta_r \rceil \rangle$.

\begin{defi}
  An $S^{-1}_2E$ proof is \emph{strictly
  $i$-normal} if all formulae contained in the proof are $i$-normal.
  The property ``$w$ is (the G\"odel number of) a strictly $i$-normal
  proof tree for $\Gamma \rightarrow \Delta$'' is
  $\Delta^b_1$-definable.  We write $i\Prf(w, \lceil \Gamma
  \rightarrow \Delta \rceil)$ for the $\Delta^b_1$ formula that
  defines this property.
\end{defi}

\begin{prop}\label{prop:soundness}
  Let $\Gamma\ENTAIL\Delta$ be a sequent comprised entirely of
  $i$-normal formulae, and let $u,w\in\mathbb{N}$ such that
  $i\Prf(w,\lceil\Gamma\ENTAIL\Delta\rceil)$ holds, $w\le u$, and the
  binary representation of $u$ is of the form $11\cdots 1$, that is,
  all the bits are $1$.  Then for every node $r$ of $w$, the following
  holds (where $\rho$ denotes an environment as well as its G\"odel
  number and $\Gamma_r \ENTAIL \Delta_r$ denotes the conclusion of
  the subproof which corresponds a node $r$).

\begin{multline}\label{eq:soundness}
  \forall \rho \leq \BdEnv(\lceil \Gamma_r \ENTAIL \Delta_r \rceil, u)
  \Bigl[ \Env(\rho, \lceil \Gamma_r \ENTAIL \Delta_r \rceil, u)
  \supset \\ \forall u' \leq u \circleddash r \{[\forall A \in
  \Gamma_r, \ T_i(u', \lceil A \rceil, \rho)] \IMPLY [\exists B \in
  \Delta_r, \ T_i(u' \oplus r, \lceil B \rceil,
  \rho)]\}\Bigr]\end{multline} Furthermore, this is derivable in
  $S^{i+2}_2$.
\end{prop}

\proof We prove the proposition by tree induction on $r$.  Since the
formula in \eqref{eq:soundness} is $\Pi^b_{i+2}$, our proof can be
carried out in $S^{i+2}_2$.

  Let $\rho\leq \BdEnv(\lceil \Gamma_r\ENTAIL\Delta_r\rceil,u)$,
  assume that $\Env(\rho,\lceil\Gamma_r\ENTAIL\Delta_r
  \rceil,u)$, and let $u'\leq u\ominus r$.  Note that $u \geq
  r$, since $u \geq w \geq r$.  Therefore, $u' \oplus r \leq u \ominus
  r \oplus r \leq u$, since all the bits in the binary representation
  of $u$ are 1. We use this fact throughout this proof.

We prove that if $\forall A \in \Gamma_r, \ T_i(u', \lceil A \rceil,
\rho)$ holds, then $\exists B \in \Delta_r, \ T_i(u' \oplus r, \lceil
B \rceil, \rho)$, by considering all possible forms for the last
inference in the derivation of $\Gamma_r \rightarrow \Delta_r$.

  \paragraph{Identity rule:}
  \begin{equation}\infer[\textup{Id}]{A \ENTAIL A}{}\end{equation}
  Assume that $T_i(u', \lceil A \rceil, \rho)$.  Then, by Lemma
  \ref{lem:T-upward}, $T_i(u' \oplus r, \lceil A \rceil, \rho)$.
  Hence, $r$ satisfies \eqref{eq:soundness}.

\paragraph{Axiom rule:}
\begin{equation}
     \infer[\textup{Ax},]{\Gamma(\vec{s}(\vec{a})) \ENTAIL
     \Delta(\vec{s}(\vec{a}))}{}
   \end{equation} where $\Gamma(\vec{s}(\vec{a})) \ENTAIL
   \Delta(\vec{s}(\vec{a}))$ is a substitution instance of an axiom.

   Since there are only finite many axioms, we use case analysis on
   the axiom which derives this substitution instance.  Assume that
   $\forall A\in \Gamma,T_i(u', \lceil A(\vec{s}(\vec{a})) \rceil,
   \rho)$.  Let $\Gamma(\vec{b}) \rightarrow \Delta(\vec{b})$ be
   the axiom into which the substitution was
     made.  By the assumption on $\mathcal{A}$, this
     axiom satisfies the boundedness conditions (Definition
   \ref{defn:boundedness}). Moreover, its
   standard interpretation (given in Definition
   \ref{defn:interpretation}) is derivable in $S^{i+2}_2$.

   By the first boundedness condition, all the formulae in $\Gamma$
   and $\Delta$ are basic. Let $\vec{b}=b_1,\dots,b_l$ and
   $\vec{s}(\vec{a})=s_1(\vec{a}),\dots,s_l(\vec{a})$, where
   $s_k(\vec{a})$ is the term that was substituted for the variable
   $b_k$ in the application of the axiom rule ($k=1,\dots,l$).  By the
   second boundedness condition, $b_k$ occurs in $\Gamma$, so by Lemma
   \ref{lem:T-sub} $\exists d_k \leq u'$ such that $v(\lceil
   s_k(\vec{a}) \rceil, \rho) \downarrow_{u'} d_k$ ($k =
   1, \dots, l$), hence $\forall A\in\Gamma, T_i(u', \lceil
   A(\vec{b}) \rceil, \rho[\vec{b} \mapsto \vec{d}])$.

   Let $\vec{t}_\Gamma(\vec{b})$ be the subterms of the terms that
   occur in $\Gamma(\vec{b})$, and let $\vec{t}_\Delta(\vec{b})$ be
   the subterms of the terms that occur in $\Delta(\vec{b})$.  Since
   all formulae occur in $\Gamma$ and $\Delta$ are basic, $\vec{b}$
   are all variables contained $\vec{t}_\Gamma(\vec{b})$ and
   $\Gamma(\vec{b})$.  Since the function symbol of $S^i_2E$ is
   definable in $S^1_2$, we can view the terms in
   $\vec{t}_\Gamma(\vec{b})$ and $\vec{t}_\Delta(\vec{b})$ as terms of
   $S^1_2$.  By the third boundedness condition, the
   relation \begin{equation} \max\{\vec{t}_\Delta(\vec{b})\} \leq \alpha
     \cdot\max\{\vec{t}_\Gamma(\vec{b}) \}
   \end{equation} is provable in $S^1_2$.

   Since $\forall A \in \Gamma, T_i(u', \lceil A(\vec{b}) \rceil,
   \rho[\vec{b} \mapsto \vec{d}])$, we have
   $\max\{\vec{t}_\Gamma(\vec{d})\} \leq u'$.  By Lemma
   \ref{lem:T0-refl}, we have that, for every $A$ in $\Gamma$,
   $A(\vec{d})$ is true (in the meta-language).  Since
   $\Gamma(\vec{d}) \ENTAIL \Delta(\vec{d})$ holds (in the
   meta-language), there is some $B$ in $\Delta$ such that
   $B(\vec{d})$ is true (in the meta-language). Since we can take
   $\alpha$ to be 4, we have $\max\{\vec{t}_\Delta(\vec{d}) \} \leq 4
   \cdot u' \leq u' \oplus r$.  Let $\vec{c} = FV(B(\vec{b}))$.  Then
   $T_i(u' \oplus r, \lceil B(\vec{c}) \rceil, \rho[\vec{c} \mapsto
   \vec{d}]$ holds by Lemma \ref{lem:T}.  By Lemma \ref{lem:v-upward}
   and the fact that $v(\lceil s_k(\vec{a}) \rceil, \rho)
   \downarrow_{u'} d_k$ ($k = 1, \dots, l$), we obtain $v(\lceil
   s_k(\vec{a}) \rceil, \rho) \downarrow_{u' \oplus r} d_k$. Using
   that result and Lemma \ref{lem:T-sub}, we have $T_i(u' \oplus r,
   \lceil B(\vec{s}(\vec{a})) \rceil,\rho)$, so we are done.

 \paragraph{Structural rules:} 
 \begin{equation} \infer{A, \Gamma \ENTAIL \Delta}{\infer*[r_1]{\Gamma
       \ENTAIL \Delta}{}} \end{equation} Assume that $T_i(u',\lceil
 A\rceil,\rho)$ and $\forall B\in\Gamma,T_i(u',\lceil
 B\rceil,\rho)$. Let $\rho_1$ be the subsequence of $\rho$ such that
 $\Env(\rho_1,\lceil\Gamma_{r_1}\ENTAIL\Delta_{r_1}\rceil,u)$. Then
 $\rho_1\le\BdEnv(\lceil\Gamma_{r_1}\ENTAIL\Delta_{r_1}\rceil,u)$ and
 by Lemma \ref{lem:T-rho} we have $\forall B\in\Gamma,T_i(u',\lceil
 B\rceil,\rho_1)$. By the induction hypothesis applied to $r_1$
 together with the fact that $u'\le u\ominus r\le u\ominus r_1$, there
 is some $B$ in $\Delta$ such that $T_i(u'\oplus r_1,\lceil
 B\rceil,\rho_1)$ holds, hence $T_i(u'\oplus r_1,\lceil B\rceil,\rho)$
 holds as well because of Lemma \ref{lem:T-rho}. By Lemma
 \ref{lem:T-upward} and the fact that $u'\oplus r_1\le u'\oplus r$, we
 are done.

 The proof for the other weakening rule is similar, and the proofs for
the remaining structural rules are trivial.

 \paragraph{$\neg$-rules:}
   \begin{equation}
       \infer[\LEFT\neg,]{\neg p(t_1,t_2), \Gamma \rightarrow
       \Delta}{\infer*[r_1]{\Gamma
         \rightarrow \Delta, p(t_1,t_2)}{}}
   \end{equation} where $p$ is $=$ or $\leq$.

   Assume that $T_i(u',\lceil \neg p(t_1,t_2) \rceil,\rho)$ and
   $\forall A \in \Gamma, T_i(u', \lceil A \rceil, \rho)$.  Note that
   $\Env(\rho,\lceil\Gamma_{r_1}\ENTAIL\Delta_{r_1}\rceil,u)$
     holds, because the variables that occur free in
     $\Gamma_{r_1}\ENTAIL\Delta_{r_1}$ are precisely those that occur
     free in $\Gamma_r\ENTAIL\Delta_r$; hence
     $\rho\le\BdEnv(\lceil\Gamma_{r_1}\ENTAIL\Delta_{r_1}\rceil,u)$. Moreover,
   $u' \leq u \ominus r \leq u \ominus r_1$, so by the induction
   hypothesis applied to $r_1$, either $T_i(u' \oplus r_1,
   \lceil p(t_1,t_2) \rceil,\rho)$ holds or $\exists B \in \Delta,
   T_i(u' \oplus r_1, \lceil B \rceil,\rho)$.  By Lemma
   \ref{lem:T-upward} and our assumption about $\neg p(t_1,t_2)$,
   $T_i(u' \oplus r_1, \lceil \neg p(t_1,t_2) \rceil,\rho)$ holds, so
   $T_i(u'\oplus r_1,\lceil p(t_1,t_2) \rceil,\rho)$ cannot also hold.
   Hence $\exists B \in \Delta, T_i(u' \oplus r_1, \lceil B
   \rceil,\rho)$.  Since $u' \oplus r_1 \leq u' \oplus r$, by Lemma
   \ref{lem:T-upward} we are done.

\begin{equation}
  \infer[\RIGHT\neg,]{Et_1, Et_2, \Gamma \rightarrow
    \Delta, \neg p(t_1, t_2)}{\infer*[r_1]{p(t_1,t_2), \Gamma
      \rightarrow \Delta}{}}
\end{equation} where $p$ is $=$ or $\leq$.

Assume that $T_i(u', \lceil Et_1 \rceil,\rho)$, $T_i(u', \lceil Et_2
\rceil,\rho)$, and $\forall A \in \Gamma, T_i(u', \lceil A \rceil,
\rho)$.  If $\exists B \in \Delta, T_i(u' \oplus r_1, \lceil B \rceil,
\rho)$, we are done, so assume otherwise, that is, for every $B \in
\Delta$, $T_i(u' \oplus r_1, \lceil B \rceil, \rho)$ does not hold.
Note that, just as in the proof for the $\LEFT\neg$-rule,
$\Env(\rho,\lceil \Gamma_{r_1}\ENTAIL\Delta_{r_1}\rceil,u)$ holds and
$\rho\le\BdEnv(\lceil
\Gamma_{r_1}\ENTAIL\Delta_{r_1}\rceil,u)$. Moreover, $u'\leq u\ominus
r\leq u\ominus r_1$. Hence, by the induction hypothesis applied to
$r_1$, $T_i(u', \lceil p(t_1, t_2) \rceil, \rho)$ does not hold.  By
Lemma \ref{lem:EM} and our assumption about $Et_1$ and $Et_2$,
$T_i(u', \lceil \neg p(t_1, t_2) \rceil, \rho)$ does hold.  Hence by
Lemma \ref{lem:T-upward} we are done.

 \paragraph{$\AND$-rules:}
\begin{equation}
   \infer[\LEFT\AND_1]{A \AND B, \Gamma
     \rightarrow \Delta}{\infer*[r_1]{A, \Gamma
       \rightarrow \Delta}{}}
 \end{equation} Assume that $T_i(u', \lceil A\AND B
 \rceil, \rho)$ and $\forall C \in \Gamma, T_i(u', \lceil C \rceil,
 \rho)$. Note that since $A \AND B$ is an $i$-normal formula, by
 Definition \ref{defn:normal} it is
 quantifier free, hence $T_i(u', \lceil A \AND B
 \rceil, \rho) \leftrightarrow T_{-1}(u', \lceil A \AND B \rceil,
 \rho)$.  By Lemma \ref{lem:T_{-1}}, we have $T_i(u', \lceil A \rceil,
 \rho)$. Let $\rho_1$ be the subsequence of $\rho$
 such that $\Env(\rho_1,\lceil
 \Gamma_{r_1}\ENTAIL\Delta_{r_1}\rceil,u)$. Then we have $T_i(u',
 \lceil A \rceil, \rho_1)$ and $\forall C \in \Gamma, T_i(u', \lceil C
 \rceil, \rho_1)$ by Lemma \ref{lem:T-rho}; in addition, $\rho_1\le \BdEnv(\lceil
 \Gamma_{r_1}\ENTAIL\Delta_{r_1}\rceil,u)$. By the induction
 hypothesis applied to $r_1$ together with the fact that $u'\leq
 u\ominus r\leq u\ominus r_1$, there
 is some $D\in\Delta$ such that $T_i(u' \oplus r_1, \lceil D \rceil, \rho_1)$
 holds.  Then $T_i(u' \oplus r_1, \lceil D \rceil, \rho)$ holds
 as well by Lemma \ref{lem:T-rho}.  Since $u'
 \oplus r_1 \leq u' \oplus r$, by Lemma \ref{lem:T-upward} we are
 done.

 The proof for the $\LEFT\AND_2$-rule is similar.

\begin{equation}
  \infer[\RIGHT\AND]{\Gamma \rightarrow \Delta, A \AND
    B}{\infer*[r_1]{\Gamma \rightarrow \Delta, A }{} \quad
    \infer*[r_2]{\Gamma \rightarrow \Delta, B}{}}
\end{equation} Assume that $\forall C \in \Gamma, T_i(u', \lceil C
\rceil, \rho)$. Let $\rho_1$ be the subsequence of $\rho$ such
that
$\Env(\rho_1,\lceil\Gamma_{r_1}\ENTAIL\Delta_{r_1}\rceil,u)$. Then we
have $\forall C \in \Gamma, T_i(u', \lceil C \rceil, \rho_1)$, and
$\rho_1\le\BdEnv(\lceil\Gamma_{r_1}\ENTAIL\Delta_{r_1}\rceil,u)$. By
the induction hypothesis applied to $r_1$ together with the fact that
$u'\leq u\ominus r\leq u\ominus r_1$, either $\exists D \in
\Delta, T_i(u' \oplus r_1, \lceil D \rceil, \rho_1)$ or $T_i(u' \oplus r_1, \lceil A
\rceil, \rho_1)$.

Similarly, let $\rho_2$ be the subsequence of $\rho$ such that
$\Env(\rho_2,\lceil \Gamma_{r_2}\ENTAIL\Delta_{r_2}\rceil,u)$. Then we
have $\forall C \in \Gamma, T_i(u', \lceil C \rceil, \rho_2)$, and
$\rho_2\le\BdEnv(\lceil\Gamma_{r_2}\ENTAIL\Delta_{r_2}\rceil,u)$. By
the induction hypothesis applied to $r_2$ together with the fact that
$u'\leq u\ominus r\leq u\ominus r_2$, either $\exists D \in \Delta,
T_i(u' \oplus r_2, \lceil D \rceil, \rho_2)$ or $T_i(u' \oplus r_2,
\lceil B \rceil, \rho_2)$.

If there exist $j\in\{1,2\}$ and $D\in\Delta$ such that $T_i(u' \oplus
r_j, \lceil D \rceil, \rho_j)$, then $T_i(u' \oplus r_j, \lceil D
\rceil, \rho)$ holds by Lemma \ref{lem:T-rho}. Thus we are done by
Lemma \ref{lem:T-upward} and the fact that $u'\oplus r_j\leq u'\oplus
r$, so assume otherwise.

Then both $T_i(u' \oplus r_1, \lceil A \rceil, \rho_1)$ and $T_i(u'
\oplus r_2, \lceil B \rceil, \rho_2)$ hold. Thus we have both $T_i(u'
\oplus r_1, \lceil A \rceil, \rho)$ and $T_i(u' \oplus r_2, \lceil B
\rceil, \rho)$, by Lemma \ref{lem:T-rho}.  Since $u' \oplus r_1, u'
\oplus r_2 \leq u' \oplus r$, both $T_i(u' \oplus r, \lceil A \rceil,
\rho)$ and $T_i(u' \oplus r, \lceil B \rceil, \rho)$ hold by Lemma
\ref{lem:T-upward}.  As noted in the proof for the $\LEFT\AND_1$-rule,
the formula $A\AND B$ is quantifier free. Thus by Lemma
\ref{lem:T_{-1}} and the definition of $T_i$, we have $T_i(u' \oplus
r, \lceil A \AND B \rceil, \rho)$.

\paragraph{$\OR$-rules:}  The proofs for $\OR$-rules are similar to
the proofs for $\AND$-rules.

\paragraph{Bounded $\forall$-rules:}
\begin{equation}
  \infer[\LEFTb\forall,]{t \leq s, \forall x \leq s.  A(x), \Gamma
  \rightarrow
    \Delta}{\infer*[r_1]{A(t), \Gamma \rightarrow \Delta}{}}
\end{equation}where the variable $x$ does
not occur in the term $s$.

Assume that $T_i(u',\lceil t\leq s\rceil,\rho)$, $T_i(u',\lceil\forall
x\leq s.A(x)\rceil,\rho)$, and $\forall B \in \Gamma,T_i(u', \lceil B
\rceil, \rho)$.  Since $T_i(u', \lceil t \leq s
\rceil, \rho)$ holds, there are $c_0, d$ such that $v(\lceil t
\rceil,\rho) \downarrow_{u'} c_0$, $v(\lceil s \rceil,\rho)
\downarrow_{u'} d$, and $c_0 \leq d$.  By Lemma \ref{lem:T} and the
fact that $v(\lceil s \rceil,\rho) \downarrow_{u'} d$, we have
$\forall c \leq d, T_i(u', \lceil A(a) \rceil, \rho[a
\mapsto c])$.  In particular, $T_i(u', \lceil A(a) \rceil, \rho[a \mapsto
c_0])$, since $c_0 \leq d$.  If $a$ occurs free in $A(a)$, we can
apply Lemma \ref{lem:T-sub} to $A(a)$ and obtain $T_i(u',
\lceil A(t) \rceil, \rho)$.  The conclusion is obvious
if $a$ does not occur free in $A$.

Note that $\rho$ is an environment for
$\Gamma_{r_1}\ENTAIL\Delta_{r_1}$, since every variable that occurs
free in $\Gamma_{r_1}\ENTAIL\Delta_{r_1}$ also occurs free in
$\Gamma_{r}\ENTAIL\Delta_r$, so let $\rho_1$ be the subsequence of
$\rho$ such that $\Env(\rho_1,\lceil
\Gamma_{r_1}\ENTAIL\Delta_{r_1}\rceil,u)$. Then we have $T_i(u',
\lceil A(t) \rceil, \rho_1)$ and $\forall B \in \Gamma,T_i(u', \lceil
B \rceil, \rho_1)$; in addition, $\rho_1\le\BdEnv(\lceil
\Gamma_{r_1}\ENTAIL\Delta_{r_1}\rceil,u)$. Thus there is
some $C$ in $\Delta$ such that $T_i(u' \oplus r_1, \lceil C \rceil,
\rho_1)$ holds, by the induction hypothesis applied to $r_1$ together
with the fact that $u'\leq u\ominus r\leq u\ominus r_1$.  Hence
we have $T_i(u' \oplus r_1, \lceil C \rceil, \rho)$, because $\rho_1$
is a subsequence of $\rho$. Since $u' \oplus r_1 \leq u' \oplus r$,
by Lemma \ref{lem:T-upward} we are done.

\begin{equation}
  \infer[\RIGHTb\forall,]{Et, \Gamma \rightarrow \Delta, \forall x
  \leq t. A(x)}{
    \infer*[r_1]{a \leq t, \Gamma \rightarrow \Delta, A(a)}{}}
\end{equation} where neither the variable $a$ nor the
variable $x$ occurs in the term
$t$, and $a$ does not occur free in
$\Gamma\ENTAIL\Delta$.

Assume that $T_i(u', \lceil Et \rceil, \rho)$ and $\forall B \in
\Gamma, T_i(u', \lceil B \rceil, \rho)$.  Since $T_i(u', \lceil Et
\rceil, \rho)$ holds, there is some $c\le u'$ such that $v(\lceil
t\rceil,\rho)\downarrow_{u'}c$.  Let $d$ be any natural number such
that $d \leq c$.  Then $T_i(u', \lceil a \leq t \rceil, \rho[a \mapsto
d])$ holds by Lemma \ref{lem:T0-refl} and the fact that $T_i(u',
\lceil a \leq t \rceil, \rho[a \mapsto d])\leftrightarrow T_{-1}(u',
\lceil a \leq t \rceil, \rho[a \mapsto d])$. Furthermore, since $a$
does not occur free in $\Gamma$, we have $\forall B \in
  \Gamma, T_i(u', \lceil B \rceil, \rho[a \mapsto d])$.

  Note that $\rho[a\mapsto d]$ is an environment for
  $\Gamma_{r_1}\ENTAIL\Delta_{r_1}$, since every variable other than
  $a$ that occurs free in $\Gamma_{r_1}\ENTAIL\Delta_{r_1}$ also
  occurs free in $\Gamma_r\ENTAIL\Delta_r$. Moreover, $\rho[a\mapsto
  d](y)\le u$ for every variable $y$ that occurs free in
  $\Gamma_{r_1}\ENTAIL\Delta_{r_1}$, since $d \leq c \leq u' \leq u$.
  Since $\Env(\rho,\lceil \Gamma_{r}\ENTAIL\Delta_{r}\rceil,u)$,
  $\Env(\rho[a \mapsto d],\lceil
  \Gamma_{r_1}\ENTAIL\Delta_{r_1}\rceil,u)$ holds.  By the induction
  hypothesis applied to $r_1$ together with the fact that $u'\leq
  u\ominus r\leq u\ominus r_1$, either $\exists C \in \Delta, T_i(u'
  \oplus r_1, \lceil C \rceil, \rho[a \mapsto d])$ or $T_i(u'\oplus
r_1, \lceil A(a) \rceil, \rho[a \mapsto d])$.

If there is some $C$ in $\Delta$ such that $T_i(u' \oplus r_1, \lceil
C \rceil, \rho[a \mapsto d])$, then we have $T_i(u' \oplus r_1,
\lceil C \rceil, \rho)$, since $a$ does not occur free in $D$. Hence
we are done by Lemma \ref{lem:T-upward} and the fact that $u'\oplus
r_1\leq u'\oplus r$, so assume otherwise.

Then $T_i(u' \oplus r_1, \lceil A(a) \rceil, \rho[a \mapsto
d])$ holds.  Since $d$ was an arbitrary natural number less
than or equal to $c$, by Lemma \ref{lem:T} we have $T_i(u'\oplus
r_1, \lceil \forall x \leq t. A(x) \rceil, \rho)$.  Since $u'\oplus
r_1\leq u'\oplus r$, we are done by Lemma \ref{lem:T-upward}.

\paragraph{Bounded $\exists$-rules:}
\begin{equation}
  \infer[\LEFTb\exists,]{\exists x \leq t. A(x),
    \Gamma \rightarrow \Delta}{\infer*[r_1]{a \leq
      t, A(a), \Gamma \rightarrow \Delta}{}}
\end{equation} where neither the variable $a$ nor the
variable $x$ occurs in the term
$t$, and $a$ does not occur free in
$\Gamma\ENTAIL\Delta$.

Assume that $T_i(u', \lceil \exists x \leq t. A(x) \rceil, \rho)$ and
$\forall B \in \Gamma, T_i(u', \lceil B \rceil, \rho)$.
By Lemma \ref{lem:T}, there exist $c,d$ such that $d\leq c$, $v(\lceil
t \rceil, \rho) \downarrow_{u'} c$, and $T_i(u', \lceil A(a) \rceil ,
\rho[a \mapsto d])$.  Then $T_i(u',\lceil a\le
t\rceil,\rho[a\mapsto d])$ holds by Lemma \ref{lem:T0-refl} and the
fact that $T_i(u', \lceil a \leq t \rceil, \rho[a \mapsto d])\equiv
T_{-1}(u', \lceil a \leq t \rceil, \rho[a \mapsto d])$. Since $a$
does not occur free in $\Gamma$, we have $\forall B \in
\Gamma, T_i(u', \lceil B \rceil, \rho[a \mapsto d])$.
 
Note that $\rho[a\mapsto d]$ is an environment for
$\Gamma_{r_1}\ENTAIL\Delta_{r_1}$, and that $\rho[a\mapsto d](y)\le u$
for every variable $y$ that occurs free in
$\Gamma_{r_1}\ENTAIL\Delta_{r_1}$, since $d \leq c \leq u' \leq
u$. Since $\Env(\rho, \lceil \Gamma_{r}\ENTAIL\Delta_{r} \rceil,u)$,
$\Env(\rho[a \mapsto d], \lceil \Gamma_{r_1}\ENTAIL\Delta_{r_1}
\rceil,u)$.  Thus by the induction
hypothesis applied to $r_1$ together with the fact that $u'\leq
u\ominus r\leq u\ominus r_1$, there is some
$C$ in $\Delta$ such that $T_i(u' \oplus r_1, \lceil C \rceil,
\rho[a \mapsto d])$.  Since $a$ does not occur free in
$C$, we have $T_i(u' \oplus r_1, \lceil C \rceil, \rho)$.  Therefore,
we are done, by Lemma \ref{lem:T-upward} and the fact that $u' \oplus
r_1 \leq u' \oplus r$.

\begin{equation}
  \infer[\RIGHTb\exists,]{t \leq s, \Gamma \rightarrow \Delta, \exists
  x \leq
    s. A(x)}{ \infer*[r_1]{\Gamma \rightarrow \Delta, A(t)}{}}
\end{equation}where the variable $x$ does
not occur in the term $s$.

Assume that $T_i(u', \lceil t \leq s \rceil, \rho)$ and $\forall B \in
\Gamma, T_i(u', \lceil B \rceil, \rho)$. Note that $\rho$ is an
environment for $\Gamma_{r_1}\ENTAIL\Delta_{r_1}$. Let $\rho_1$ be the
subsequence of $\rho$ such that
$\Env(\rho_1,\lceil\Gamma_{r_1}\ENTAIL\Delta_{r_1}\rceil,u)$. Then we
have $\forall B \in \Gamma, T_i(u', \lceil B \rceil, \rho_1)$, and
$\rho_1\le\BdEnv(\lceil\Gamma_{r_1}\ENTAIL\Delta_{r_1}\rceil,u)$. By
the induction hypothesis applied to $r_1$ together with the fact that
$u'\leq u\ominus r\leq u\ominus r_1$, either $\exists C \in \Delta,
T_i(u' \oplus r_1, \lceil C \rceil, \rho_1)$ or $T_i(u' \oplus r_1,
\lceil A(t) \rceil, \rho_1)$.

If there is some $C$ in $\Delta$ such that $T_i(u' \oplus r_1, \lceil
C \rceil, \rho_1)$, then we have $T_i(u' \oplus r_1, \lceil C \rceil,
\rho)$ by Lemma \ref{lem:T-rho}. Hence we are done by Lemma
\ref{lem:T-upward} and the fact that $u'\oplus r_1\leq u'\oplus r$, so
assume otherwise.

Then we have $T_i(u' \oplus r_1, \lceil A(t) \rceil, \rho_1)$, hence
$T_i(u' \oplus r_1, \lceil A(t) \rceil, \rho)$ holds by Lemma
\ref{lem:T-rho}. By Lemma \ref{lem:T-sub}, there exists $c$ such that
$v(\lceil t \rceil, \rho)\downarrow_{u' \oplus r_1} c$ and $T_i(u'
\oplus r_1, \lceil A(a) \rceil, \rho[a \mapsto c])$.  Since $T_i(u',
\lceil t \leq s \rceil, \rho)$ holds, there exists $d$ such that
$v(\lceil s \rceil, \rho)\downarrow_{u'} d$ and $c \leq d$. By Lemma
\ref{lem:T}, $T_i(u' \oplus r_1, \lceil \exists x \leq s. A(x) \rceil,
\rho)$ holds. Thus we are done, by Lemma \ref{lem:T-upward} and the
fact that $u'\oplus r_1\leq u'\oplus r$.

\paragraph{Cut rule:}
\begin{equation}\label{eq:Cutr1r2}
  \infer[\textup{Cut}]{\Gamma, \Pi \rightarrow \Delta, \Lambda}{
    \infer*[r_1]{\Gamma \rightarrow \Delta, A}{} \quad
    \infer*[r_2]{A, \Pi \rightarrow \Lambda}{}
  }
\end{equation} 
Assume that $\forall B \in \Gamma, \Pi, T_i(u', \lceil B \rceil,
\rho)$.  Let $\vec{a}$ be the variables that occur free in $A$ but do
not occur free in $\Gamma_r\ENTAIL\Delta_r$, and let
$\rho[\vec{a}\mapsto \vec{0}]$ be the environment that extends $\rho$
and maps every variable in $\vec{a}$ to 0 (where
$\rho[\vec{a}\mapsto\vec{0}]\equiv\rho$ if $\vec{a}$ is empty). Note
that $\rho[\vec{a}\mapsto\vec{0}]$ is an environment for
$\Gamma_{r_1}\ENTAIL\Delta_{r_1}$, and that
$\rho[\vec{a}\mapsto\vec{0}](y)\le u$ for every variable $y$ that
occurs free in $\Gamma_{r_1}\ENTAIL\Delta_{r_1}$. Let $\rho_1$ be the
subsequence of $\rho[\vec{a}\mapsto \vec{0}]$ such that
$\Env(\rho_1,\lceil\Gamma_{r_1}\ENTAIL\Delta_{r_1}\rceil,u)$. Then we
have $\forall B \in \Gamma, T_i(u', \lceil B \rceil, \rho_1)$, and
$\rho_1\le\BdEnv(\lceil\Gamma_{r_1}\ENTAIL\Delta_{r_1}\rceil,u)$. By
the induction hypothesis applied to $r_1$ together with the fact that
$u'\leq u\ominus r\leq u\ominus r_1$, either $\exists C \in \Delta,
T_i(u' \oplus r_1, \lceil C \rceil, \rho_1)$ or $T_i(u' \oplus r_1,
\lceil A \rceil, \rho_1)$.

If there is some $C$ in $\Delta$ such that $T_i(u' \oplus r_1, \lceil
C \rceil, \rho_1)$, then we have $T_i(u' \oplus r_1, \lceil C \rceil,
\rho[\vec{a}\mapsto \vec{0}])$ because $\rho_1$ is a subsequence of
$\rho[\vec{a}\mapsto \vec{0}]$ and by Lemma
\ref{lem:T-rho}. Furthermore, $T_i(u' \oplus r_1, \lceil C \rceil,
\rho)$ holds by Lemma \ref{lem:T-rho}, since none of the variables in
$\vec{a}$ occurs free in $C$. Thus $T_i(u' \oplus r, \lceil C \rceil,
\rho)$ by Lemma \ref{lem:T-upward} and the fact that $u'\oplus r_1\leq
u'\oplus r$. Hence we are done, so assume otherwise.

Then we have $T_i(u' \oplus r_1, \lceil A \rceil, \rho_1)$, so $T_i(u'
\oplus r_1, \lceil A \rceil, \rho[\vec{a}\mapsto \vec{0}])$ holds by
Lemma \ref{lem:T-rho}, because $\rho_1$ is a subsequence of
$\rho[\vec{a}\mapsto \vec{0}]$. By our assumption about $\Pi$, we have
$\forall B\in\Pi,T_i(u',\lceil B\rceil,\rho[\vec{a}\mapsto\vec{0}])$
by Lemma \ref{lem:T-rho}, because $\rho$ is a subsequence of
$\rho[\vec{a}\mapsto\vec{0}]$. By Lemma \ref{lem:T-upward}, we have
$\forall B \in \Pi, T_i(u' \oplus r_1, \lceil B \rceil,
\rho[\vec{a}\mapsto\vec{0}])$.
  
Note that $\rho[\vec{a}\mapsto \vec{0}]$ is an environment for
$\Gamma_{r_2}\ENTAIL\Delta_{r_2}$, and that $\rho[\vec{a}\mapsto
\vec{0}](y)\le u$ for every variable $y$ that occurs free in
$\Gamma_{r_2}\ENTAIL\Delta_{r_2}$. Let $\rho_2$ be the subsequence of
$\rho[\vec{a}\mapsto \vec{0}]$ such that
$\Env(\rho_2,\lceil\Gamma_{r_2}\ENTAIL\Delta_{r_2}\rceil,u)$. Then we
have $T_i(u' \oplus r_1, \lceil A \rceil, \rho_2)$ and $\forall B \in
\Pi, T_i(u' \oplus r_1, \lceil B \rceil, \rho_2)$; in addition,
$\rho_2\le\BdEnv(\lceil\Gamma_{r_2}\ENTAIL\Delta_{r_2}\rceil,u)$.

Our choice of G\"odel numbering, together with the fact that $r_1$ and
$r_2$ are G\"odel numbers of nonempty subproofs of
$\Gamma_r\ENTAIL\Delta_r$, ensures that $|r_1\oplus r_2|<|r|$.  Since
$u'\le u\ominus r$, we have $|u' \oplus r_1| \leq |u \ominus r \oplus
r_1| <|u \ominus (r_1\oplus r_2) \oplus r_1|= |u \ominus r_2|$, hence
$u'\oplus r_1<u\ominus r_2$.

By the induction hypothesis applied to $r_2$ together with the fact
that $u'\oplus r_1<u\ominus r_2$, there is some $D$ in $\Lambda$ such
that $T_i(u' \oplus r_1\oplus r_2, \lceil D \rceil, \rho_2)$. Then we
have $T_i(u' \oplus r_1\oplus r_2, \lceil D \rceil,
\rho[\vec{a}\mapsto \vec{0}])$ by Lemma \ref{lem:T-rho}, because
$\rho_2$ is a subsequence of $\rho[\vec{a}\mapsto
\vec{0}]$. Furthermore, $T_i(u' \oplus r_1\oplus r_2, \lceil D \rceil,
\rho)$ holds by Lemma \ref{lem:T-rho}, because none of the variables
in $\vec{a}$ occurs free in $D$.  Since $|r_1\oplus r_2|<|r|$, we
  have $u' \oplus r_1 \oplus r_2 < u' \oplus r$, so $T_i(u'
  \oplus r, \lceil D \rceil, \rho)$ by Lemma \ref{lem:T-upward}.
Hence we are done.
\qed

\begin{thm}
  Let $i\Con \equiv \forall w. \neg i\Prf(w, \lceil \rightarrow
  \rceil)$, which states that there is no strictly $i$-normal proof of
  the empty sequent $\ENTAIL$.  Then
\begin{equation}
  S^{i+2}_2 \vdash i\Con
\end{equation}
\end{thm}

\proof
  We informally argue inside of $S^{i+2}_2$.  Assume that $i\Prf(w,
  \lceil \ENTAIL \rceil)$ holds for some $w$.  Let $u$ be as in the
  statement of Proposition \ref{prop:soundness}, let $\rho$ be the
  empty environment, and let $r$ be the root of $w$.  Then we obtain
  $[\forall A \in \Gamma_r, \ T_i(u', \lceil A \rceil, \rho)] \IMPLY
  [\exists B \in \Delta_r, \ T_i(u' \oplus r, \lceil B \rceil,
  \rho)]$.  However, both $\Gamma_r$ and $\Delta_r$ are
  empty.  Therefore, we obtain $[\forall A
  \in \emptyset, T_i(u', \lceil A \rceil, \rho)] \IMPLY [ \exists B
  \in \emptyset, T_i(u' \oplus r, \lceil B \rceil, \rho)]$.  Since
  there is no $A \in \emptyset$, the premise is true.  But since there
  is no $B \in \emptyset$, the conclusion cannot be true.  This is a
  contradiction.  Therefore, the formula $\forall w. \neg
  i\Prf(w, \lceil \rightarrow \rceil)$ holds.
\qed

\section{Bootstrapping Theorem of $S^i_2E$}
\label{sec:bootstrapping}

In this section, we establish the correspondence between $S^i_2E$ and
$S^i_2$.  We show that $S^i_2E$ has essentially the same strength as
$S^i_2$ if $i \geq 1$.  The theorem which establishes the
correspondence is called the Bootstrapping Theorem (Theorem
\ref{thm:bootstrap}), following Buss' use of the term
``bootstrapping'' in \cite{Buss:Book}, since we bootstrap from the
restricted set of axioms of $S^i_2E$ to the full power of $S^i_2$.

We present a proof of the theorem in four ``phases'' of
bootstrapping. In the first phase, we show that all the functions of
$S_2E$ are provably total.  Each of the remaining phases applies to a
particular class of inferences of $S^i_2$, and we show that all the
inferences covered in each phase are \emph{admissible in $S^i_2E$} (if
properly translated from $S^i_2$ to $S^i_2E$), that is, that if all
the premises of an inference covered in a given phase are provable in
$S^i_2E$, then the conclusion of that inference is also provable in
$S^i_2E$ (Definition \ref{defn:admissible}). The Bootstrapping Theorem
(Theorem \ref{thm:bootstrap}) then follows from the fact that every
inference of $S^i_2$ is treated in some phase of the bootstrapping.
Even the axioms are included in this, since an axiom is just a rule of
inference with no premise.

\subsection{Translation of theorems of
  $S^i_2$}\label{subsec:translation}

In this subsection, we introduce a translation of
$S^i_2$ formulae to the language of $S^i_2E$ and state the
Bootstrapping Theorem (Theorem \ref{thm:bootstrap}).

\begin{defi}\label{defn:embed}
  The formulae of $S^i_2$ are translated into formulae of $S^i_2E$ by
  replacing every formula of the form $A \IMPLY B$ with one of the
  form $\neg A \OR B$, and using De Morgan duality to replace every
  formula of the form $\neg A$ with a logically equivalent formula in
  which every subformula prefaced with the negation symbol ``$\neg$"
  is of the form $t_1 = t_2$ or $t_1 \leq t_2$.  We call this
  translation the \emph{$*$-translation} and denote the $*$-translation
of $A$ by $A^*$.  Formally, the $*$-translation is defined as follows.
  \begin{enumerate}[(1)]
  \item $(p(t_1,t_2))^*\equiv p(t_1,t_2)$ if $p$ is $=$ or $\leq$.
  \item $(\neg p(t_1, t_2))^* \equiv \neg p(t_1, t_2)$ if $p$ is $=$
    or $\leq$.
  \item $(A \AND B)^* \equiv A^* \AND B^*$.
  \item $(A \OR B)^* \equiv A^* \OR B^*$.
  \item $(\neg A)^* \equiv (\overline{A})^*$, where $\overline{A}$ is
    the De Morgan dual of $A$.
  \item $(A \IMPLY B)^* \equiv (\overline{A})^* \OR B^*$. % ??
  \item $(\forall x \leq t. A)^* \equiv \forall x \leq t. A^*$ and
    $(\exists x \leq t. A)^* \equiv \exists x \leq t. A^*$.
  \item $(\forall x. A)^* \equiv \forall x. A^*$ and
    $(\exists x. A)^* \equiv \exists x. A^*$.
  \end{enumerate}
 
 $\Gamma^*$ is the sequence of formulae which is obtained by
  applying $*$ to the formulae in the sequence $\Gamma$.

  The sequent $\Gamma \ENTAIL \Delta$ is translated to the sequent
  $(\Gamma\ENTAIL\Delta)^*\equiv E\vec{a}, \Gamma^* \ENTAIL \Delta^*$,
  where $\vec{a}$ are the variables that occur free in
  $\Gamma\ENTAIL\Delta$.
\end{defi}

The following theorem states that $S^i_2E$ proves the $*$-translations
of sequents derivable in $S^i_2$ if $i \geq 0$.
\begin{thm}[Bootstrapping Theorem]\label{thm:bootstrap}
  If $i\geq 1$ and $S^i_2$ proves a sequent $\Gamma \ENTAIL \Delta$,
  then $S^i_2E(\mathcal{F}, \mathcal{A})$ proves its $*$-translation
  $(\Gamma \ENTAIL \Delta)^*$ if $\mathcal F$ and $\mathcal{A}$
  satisfy the conditions presented in Subsections \ref{subsec:lang}
  and \ref{subsec:axiom}, respectively.
\end{thm}

The rest of this section is devoted to a proof of the Bootstrapping
Theorem.  To simplify the notation, we write $S^i_2E$ for
$S^i_2E(\mathcal{F}, \mathcal{A})$

\subsection{Bootstrapping Phase I: $S^i_2E$ proves totality of its
  functions.}\label{subsec:bootstrap_I}

In this subsection, we prove that if $i\geq 0$, all the functions of
$S^i_2E$ are provably total, that is, that $S^i_2E \vdash E\vec{a}
\ENTAIL Ef\vec{a}$ for every function symbol $f\in\mathcal{F}$.  The
proof is by induction (in the meta-language) on the definition degree
of $f$ (Definition \ref{defn:defdeg}).

\begin{prop}\label{prop:total} If $i\geq 0$, then for every
$n$-ary function symbol $f$ of $S^i_2E$, $S^i_2E$ proves
  \begin{equation}
    E\vec{a} \ENTAIL Ef\vec{a},
  \end{equation} where $\vec{a}\equiv a_1, \ldots, a_n$.
\end{prop} The
reason for specifying that $i\ge 0$ is that in the proof we apply the
$\PIND$ rule to $\Sigma^b_0$ formulae of $S^i_2E$.

It follows from this proposition that if all the variables in a term
of $S^i_2E$ converge, then the term itself converges.
\begin{cor}\label{cor:conv}
  Let $t$ be a term of $S^i_2E$.  If $a_1, \ldots, a_n$ are the
  variables that occur in $t$, then the following holds if $i\geq 0$.
  \begin{equation}
    S^i_2E \vdash Ea_1, \ldots, Ea_n \ENTAIL Et
  \end{equation}
\end{cor}

\proof[Proof of Corollary \ref{cor:conv}]
  Induction on the construction of $t$.

  The base case $t\equiv 0$: This is immediate, since $\ENTAIL E0$ is
  Axiom \eqref{eq:Ax-data-0}.

  The base case $t \equiv a_1$: This is also immediate, by Identity
  ($Ea_1\ENTAIL Ea_1$).

  Induction step $t \equiv f t_1 \cdots t_m$: We assume that the
  corollary holds for $t_1,\dots,t_m$.  Then for every $j$, we have an
  $S^i_2E$ proof of $Ea_1, \ldots, Ea_n \ENTAIL Et_j$ by the induction
  hypothesis (together with Weakening if at least one of the variables
  $a_1,\dots,a_n$ does not occur in $t_j$).  Let $b_1,\dots,b_m$ be
  variables.  By Proposition \ref{prop:total}, we have $Eb_1, \dots,
  Eb_m \ENTAIL Efb_1 \cdots b_m$.  By the Substitution Lemma (Lemma
  \ref{lem:subst}), we have $Et_1, \dots, Et_m \ENTAIL Eft_1 \cdots
  t_m$.  Applying Cut $m$ times (once for each $j$), followed by
  Contraction every time but the first, we obtain $Ea_1, \ldots, Ea_n
  \ENTAIL Eft_1 \dots t_m$.
\qed

The rest of this subsection is devoted to a proof of Proposition
\ref{prop:total}.
\proof[Proof of Proposition \ref{prop:total}]
  The proof is by induction on $d(f)$, the definition degree
  of the function $f$.
   
  Base case: $d(f)=0$.  If $d(f)=0$, then $f$ is either the constant
  function $0^n$, a projection function $\textup{proj}^n_k$, one of
  the binary successor functions $s_0,s_1$.
  If $f$ is $0^n$, then
  \begin{equation}
    E\vec{a}, E0 \ENTAIL 0^n(a_1, \dots, a_n) = 0
  \end{equation} is an axiom (See Definition \ref{defn:defnAx}).  By
  $\ENTAIL E0$ (Axiom \eqref{eq:Ax-data-0}) together with Cut, we have
  $E\vec{a} \ENTAIL 0^n(a_1, \dots, a_n) = 0$.  Using a substitution
  instance of Axiom \eqref{eq:Ax-E-p} with $p$ set to =, we can derive
  $0^n(a_1, \dots, a_n) = 0\ENTAIL E0^n(a_1, \dots, a_n)$.  Thus, we
  obtain $Ea_1, \dots, Ea_n \ENTAIL E0^n(a_1, \dots, a_n)$ by Cut.

 If $f$ is $\textup{proj}^n_k$, then
    \begin{equation} E\vec{a} \ENTAIL \textup{proj}^n_k(a_1, \dots,
    a_n) = a_k
  \end{equation} is an axiom.  Using a substitution instance of Axiom
  \eqref{eq:Ax-E-p} with Cut, we can derive 
  \begin{equation}
   Ea_1, \dots, Ea_n \ENTAIL E\textup{proj}^n_k(a_1, \dots, a_n).
  \end{equation}

If $f$ is a binary successor function $s_j$ for some
$j\in\{0,1\}$, then $Ea\ENTAIL Es_ja$ is a data axiom (Definition
\ref{defn:data-axiom}).

Induction step: $f$ is defined either by recursion or by composition.
We first consider the case of composition.

Assume that $f$ is defined from functions $g, h_1, \ldots, h_m$
by composition.  Then the defining axiom for $f$ has the following
form.
  \begin{equation}\label{eq:E-comp}
    E\vec{a}, Eg(h_1(\vec{a}), \ldots, h_m(\vec{a})) \ENTAIL
    f(\vec{a}) = g(h_1(\vec{a}), \ldots, h_m(\vec{a}))
   \end{equation}

   By the induction hypothesis, $S^i_2E$ proves
   \begin{equation}\label{eq:Eg-comp}
     Eb_1, \ldots, Eb_m \ENTAIL Eg(b_1, \ldots, b_m),
   \end{equation}where $b_1,\dots,b_m$ are variables not in
   $\{a_1,\dots,a_n\}$, and
   \begin{equation}\label{eq:Eh-comp}
     Ea_1, \ldots, Ea_n \ENTAIL Eh_j(a_1, \ldots, a_n)
   \end{equation} for $j\in\{1,\ldots,m\}$.

   Substituting $h_j(a_1, \ldots, a_n)$ for $b_j$ ($j=1,\dots,m$) in
   \eqref{eq:Eg-comp} and using the Substitution Lemma (Lemma
   \ref{lem:subst}), we obtain an $S^i_2E$
   proof of
   \begin{multline}\label{eq:Egh}
     Eh_1(a_1, \ldots, a_n), \ldots, Eh_m(a_1, \ldots, a_n) \ENTAIL\\
     Eg(h_1(a_1, \ldots, a_n), \ldots, h_m(a_1, \ldots, a_n)).
   \end{multline}

   Applying Cut to \eqref{eq:Egh} and \eqref{eq:Eh-comp} $m$ times
   (once for every $j$), followed by Contraction every time but the
   first, we have an $S^i_2E$ proof of
\begin{equation}\label{eq:Egh-x}
  Ea_1, \ldots, Ea_n \ENTAIL Eg(h_1(a_1, \ldots, a_n), \ldots,
  h_m(a_1, \ldots, a_n)).
\end{equation}

Applying Cut to \eqref{eq:Egh-x} and \eqref{eq:E-comp},
   followed by Contraction, we have
\begin{equation} Ea_1, \ldots, Ea_n \ENTAIL f(\vec{a}) =
g(h_1(\vec{a}), \ldots, h_m(\vec{a}))
\end{equation} Finally, using a substitution instance of Axiom
\eqref{eq:Ax-E-p} together with Cut, we obtain an
$S^i_2E$ proof of
\begin{equation}
  Ea_1, \ldots, Ea_n \ENTAIL Ef(a_1, \ldots, a_n).
\end{equation}

Next, we consider the case where $f$ is defined from functions $g,
h_0, h_1$ by recursion.  Then the defining
axioms for $f$ have the following forms.
\begin{align}
  E\vec{a}, Eg(\vec{a}) &\ENTAIL f(0, \vec{a}) =
  g(\vec{a}) \label{eq:E-rec-g}\\
  Ea,Eh_j(a, f(a, \vec{a}), \vec{a}) &\ENTAIL f(s_ja, \vec{a}) =
  h_j(a, f(a, \vec{a}), \vec{a}) 
   \label{eq:E-rec-h}
 \end{align}
 where $j\in\{0,1\}$.
  
  By the induction hypothesis, $S^i_2E$ proves
  \begin{equation}\label{eq:Eg-rec}
    Ea_1, \ldots, Ea_n \ENTAIL Eg(a_1,\ldots, a_n)
  \end{equation} and
\begin{equation}\label{eq:Ef-rec}
    Ea, Eb, Ea_1, \ldots, Ea_n \ENTAIL Eh_j(a, b, a_1, \ldots, a_n),
  \end{equation}where $b$ is a variable not in
  $\{a,a_1,\dots,a_n\}$ and $j\in\{0,1\}$.  Applying Cut to
  \eqref{eq:Eg-rec} and \eqref{eq:E-rec-g}, followed by Contraction,
  we obtain
\begin{equation} Ea_1,\dots,Ea_n \ENTAIL f(0, a_1,\dots,a_n) =
g(a_1,\dots,a_n)
\end{equation}
Using a substitution instance of Axiom \eqref{eq:Ax-E-p}
and Cut, we can derive
  \begin{equation}\label{eq:Ef0}
    Ea_1, \ldots, Ea_n \ENTAIL Ef(0, a_1, \ldots, a_n)
  \end{equation}

  Substituting $f(a, a_1, \ldots, a_n)$ for $b$ in \eqref{eq:Ef-rec}
  and using the Substitution Lemma (Lemma \ref{lem:subst}), we
  obtain $S^i_2E$ proofs of
  \begin{multline}
   Ea, Ef(a, a_1, \ldots, a_n), Ea_1, \ldots, Ea_n \ENTAIL\\ Eh_j(a,
   f(a, a_1, \ldots, a_n), a_1,
    \ldots, a_n)
  \end{multline}for $j\in\{0,1\}$.

Applying Cut to this result and \eqref{eq:E-rec-h}, together with totality of $\Cond$
$\ominus$, and $k$, we obtain
\begin{equation} Ea, Ef(a, \vec{a}), E\vec{a}\ENTAIL f(s_ja, \vec{a})
  = h_j(a, f(a, \vec{a}), \vec{a})
\end{equation}for $j\in\{0,1\}$. Using a
substitution instance of Axiom \eqref{eq:Ax-E-p} together with Cut, we
have $S^i_2E$ proofs of
  \begin{align}
    Ea, Ef(a, a_1, \ldots, a_n), Ea_1, \ldots, Ea_n &\ENTAIL Ef(s_0a,
    a_1, \ldots, a_n) \label{eq:Efs0}\\ Ea, Ef(a, a_1, \ldots, a_n),
    Ea_1, \ldots, Ea_n &\ENTAIL Ef(s_1a, a_1, \ldots,
    a_n). \label{eq:Efs1}
  \end{align}

  Applying the $\Sigma^b_0\text{-}\PIND\text{-}E$ rule
  \eqref{eq:PIND-E} to \eqref{eq:Ef0}, \eqref{eq:Efs0}, and
  \eqref{eq:Efs1}, and setting $t$ to
  $a$, we have an $S^i_2E$ proof of
  \begin{equation}
    Ea, Ea_1, \ldots, Ea_n \ENTAIL Ef(a, a_1, \ldots, a_n).
  \end{equation}

This completes the induction step.

\qed

\subsection{Bootstrapping Phase II : $S^i_2E$ proves $*$-translations
  of axioms of $S^i_2$}\label{sec:bootstrap_II}

In Bootstrapping Phase II, we prove the
$*$-translations of axioms of $S^i_2$ in
$S^i_2E$.  There are two kinds of axioms: equality axioms and BASIC
axioms.

\begin{prop}\label{prop:equality}
  The $*$-translations of the equality axioms of $S^i_2$ are provable
  in $S^i_2E$.
\end{prop}

\proof
  First, we consider the equality axiom $\vec{a} = \vec{b} \ENTAIL
  f(\vec{a}) = f(\vec{b})$.  The $*$-translation of this is $E\vec{a},
  E\vec{b}, \vec{a} = \vec{b} \ENTAIL f(\vec{a}) = f(\vec{b})$.  By
  Proposition \ref{prop:total}, $S^i_2E$ proves $E\vec{a} \ENTAIL
  Ef(\vec{a})$.  Therefore, we can derive it from Axiom
  \eqref{eq:Ax-eq-fun}.

  The proofs for the other equality axioms are straightforward.
\qed

Next we prove that the $*$-translations of the BASIC axioms of $S^i_2$
are provable in $S^i_2E$.
\begin{prop}\label{prop:BASIC}
  Assume that $A$ is a BASIC axiom.  Then $(\ENTAIL A)^*$ (the
  $*$-translation of $\ENTAIL A$) is derivable in $S^i_2E$.
\end{prop}

\proof
  The BASIC axioms can be derived from the auxiliary axioms by
  Corollary \ref{cor:conv} and propositional inferences.  For example,
  we consider the BASIC axiom $|a| = |b| \IMPLY a\#c = b \# c$.  We
  have the corresponding auxiliary axiom $|a| = |b|, Ea\#c, Eb\#c
  \ENTAIL a\#c = b \# c$.  By Corollary \ref{cor:conv}, we have $Ea,
  Ec \ENTAIL Ea\#c$ and $Eb, Ec \ENTAIL Eb\#c$.  Therefore, we have
  $|a| = |b|, Ea, Eb, Ec \ENTAIL a\#c = b \# c$ by Cut and
  Contraction.  By propositional inference, $Ea, Eb, Ec, E|a|, E|b|
  \ENTAIL |a| \not= |b| \OR a\#c = b \# c$.  Again, by Corollary
  \ref{cor:conv} and Cut, we have $Ea, Eb, Ec \ENTAIL |a| \not= |b|
  \OR a\#c = b \# c$, which is the $*$-translation of $\ENTAIL |a| =
  |b| \IMPLY a\#c = b \# c$.
\qed

\subsection{Bootstrapping Phase III :$*$-translations of predicate
  logic are admissible in $S^i_2E$}
\label{subsec:bootstrap_III}

In Bootstrapping Phase III, we prove that the $*$-translations of the
inferences of predicate logic are admissible in $S^i_2E$.

\begin{defi}\label{defn:admissible}
  The inference \begin{equation}
    \infer[]{\Gamma \ENTAIL \Delta}{
      \Gamma_1 \ENTAIL \Delta_1\quad \cdots\quad \Gamma_n \ENTAIL
    \Delta_n }
  \end{equation} is \emph{admissible in $S^i_2E$} if $\Gamma \ENTAIL
  \Delta$ is provable in $S^i_2E$ whenever $\Gamma_1 \ENTAIL \Delta_1,
  \dots, \Gamma_n \ENTAIL \Delta_n$ are provable in $S^i_2E$.
\end{defi}

\begin{prop}\label{prop:pred}
  If
  \begin{equation}
    \infer{\Gamma \ENTAIL \Delta}{
      \Gamma_1 \ENTAIL \Delta_1 \quad \cdots \quad \Gamma_n \ENTAIL
    \Delta_n }
  \end{equation} is an inference of predicate logic, then the
  inference\begin{equation}
    \infer{(\Gamma \ENTAIL \Delta)^*}{
      (\Gamma_1 \ENTAIL \Delta_1)^* \quad \cdots \quad (\Gamma_n
      \ENTAIL \Delta_n)^*
    }
  \end{equation} is admissible in $S^i_2E$.
\end{prop}

We prove the proposition by considering the various rules of inference
of $S^i_2$.  We begin by providing detailed proofs
for the $\LEFT\neg$-rule and the $\RIGHT\neg$-rule, in Lemma
\ref{lem:Lneg} and Lemma \ref{lem:Rneg}, respectively. Then we proceed
to proofs for the $\IMPLY$-rules and the quantifier rules.  The proofs
for the other rules are trivial.

\begin{lem}\label{lem:Lneg}
  \begin{equation}
    \infer{E\vec{a}, (\neg A)^*, \Gamma^* \ENTAIL \Delta^*}{
      E\vec{a}, \Gamma^* \ENTAIL \Delta^*, A^* }
  \end{equation} is admissible in $S^i_2E$, where $\vec{a}$ are the variables that
  occur free in the sequent $\Gamma\ENTAIL\Delta,A$.
\end{lem}

\proof By induction on $A$.

If $A$ is atomic, the inference given in the statement of the Lemma is
an instance of the $\LEFT\neg$-rule of $S^i_2E$.  In the induction
step, we assume that $E\vec{a}, \Gamma^* \ENTAIL \Delta^*, A^*$ holds
and show that $E\vec{a}, (\neg A)^*, \Gamma^* \ENTAIL \Delta^*$ holds.
The proof depends on the form of $A$.
\paragraph{$A \equiv \neg A_1$:} By Identity and Weakening, we have
$E\vec{b},(A_1)^*\ENTAIL (A_1)^*$, where $\vec{b}$ are the variables
that occur free in $A_1$. By the induction hypothesis applied to
$A_1$, we can prove $E\vec{b}, (A_1)^*, (\neg A_1)^* \ENTAIL $ in
$S^i_2E$. By assumption, $E\vec{a}, \Gamma^* \ENTAIL \Delta^*, (\neg
A_1)^* $ is derivable. Thus $E\vec{a}, (A_1)^*, \Gamma^* \ENTAIL
\Delta^*$ is derivable by Cut and Contraction.  By the definition of
$*$, we are done. (Note that, by De Morgan duality, $(\neg \neg
A_1)^*\equiv (A_1)^*$.)
  
\paragraph{$A \equiv A_1 \AND A_2$:} Since $(A_1 \AND A_2)^* \ENTAIL
(A_1)^*$ and $(A_1 \AND A_2)^* \ENTAIL (A_2)^*$ are derivable in
$S^i_2E$ by purely propositional reasoning, $E\vec{a}, \Gamma^*
\ENTAIL \Delta^*, (A_1)^*$ and $E\vec{a}, \Gamma^* \ENTAIL \Delta^*,
(A_2)^*$ are derivable from the assumption that $E\vec{a}, \Gamma^*
\ENTAIL \Delta, (A_1)^* \AND (A_2)^*$ is provable.  By the induction
hypothesis applied to $A_1$ and $A_2$, $E\vec{a}, (\neg A_1)^*,
\Gamma^* \ENTAIL \Delta^*$ and $E\vec{a}, (\neg A_2)^*, \Gamma^*
\ENTAIL \Delta^*$ are derivable.  Thus by the $\LEFT\OR$-rule of
$S^i_2E$, we can derive $E\vec{a}, (\neg A_1)^* \OR (\neg A_2)^*,
\Gamma^* \ENTAIL \Delta^*$.  By the definition of $*$, we are done.

\paragraph{$A \equiv A_1 \OR A_2$:} Since $(A_1)^* \OR (A_2)^* \ENTAIL
(A_1)^*, (A_2)^*$ is derivable by purely propositional reasoning in
$S^i_2E$, we have $E\vec{a}, \Gamma^* \ENTAIL \Delta^*, (A_1)^*,
(A_2)^*$ from the assumption that $E\vec{a}, \Gamma^* \ENTAIL
\Delta^*, (A_1)^* \OR (A_2)^*$.  Then $E\vec{a}, (\neg A_1)^*,
\Gamma^* \ENTAIL \Delta^*,(A_2)^*$, by the induction hypothesis
applied to $A_1$.  Again applying the induction hypothesis, this time
to $A_2$, we obtain a sequent $E\vec{a}, \neg (A_1)^*, \neg (A_2)^*, \Gamma^*
\ENTAIL \Delta^*$.  By the $\LEFT\AND_1$- and $\LEFT\AND_2$-rules of
$S^i_2E$, together with Contraction, we have $E\vec{a}, \neg (A_1)^*
\AND \neg (A_2)^*, \Gamma^* \ENTAIL \Delta^*$.  By the definition of
$*$, we are done.

\paragraph{$A \equiv A_1 \IMPLY A_2$:} Since $(\neg A_1)^* \OR (A_2)^*
\ENTAIL (\neg A_1)^*, (A_2)^*$ is derivable by purely propositional
reasoning in $S^i_2E$, we have $E\vec{a}, \Gamma^* \ENTAIL \Delta^*,
(\neg A_1)^*, (A_2)^*$ from the assumption that $E\vec{a}, \Gamma^*
\ENTAIL \Delta^*, (\neg A_1)^* \OR (A_2)^*$.  Applying the induction
hypothesis twice in succession (once to $\neg A_1$ and once to $A_2$),
we have $E\vec{a}, (A_1)^*, (\neg A_2)^*, \Gamma^* \ENTAIL \Delta^*$
(since $(\neg \neg A_1)^* \equiv (A_1)^*$). By the $\LEFT\AND_1$- and
$\LEFT\AND_2$-rules of $S^i_2E$, together with Contraction, we have
$E\vec{a}, (A_1)^* \AND (\neg A_2)^*, \Gamma^* \ENTAIL \Delta^*$.  By
definition of $*$, we are done.

\paragraph{$A \equiv \forall x \leq t. A_1(x)$:} Let $a$ be a variable
that occurs in neither $A$ nor $\Gamma\ENTAIL\Delta$. By Identity, we
have $A_1(a)^*\ENTAIL A_1(a)^*$. Thus $E\vec{b},Ea,a=a,A_1(a)^*\ENTAIL
A_1(a)^*$ holds by Weakening, where $\vec{b}$ are the variables other
than $a$ that occur free in $A_1(a)$.  Since $a$ occurs free in
$a=a,A_1(a)^*$, we can apply the induction hypothesis to
$E\vec{b},Ea,a=a,A_1(a)^*\ENTAIL A_1(a)^*$, hence we have
$E\vec{b},Ea,a=a,A_1(a)^*,(\neg A_1(a))^*\ENTAIL$.

Using $Ea\ENTAIL a=a$ (Axiom \eqref{eq:Ax-eq-refl}), we can derive
$E\vec{b},Ea,A_1(a)^*,(\neg A_1(a))^*\ENTAIL$ by Cut and
Contraction. By the $\LEFTb\forall$-rule of $S^i_2E$, we obtain
\begin{equation}
E\vec{b},Ea,a\leq t,\forall x\leq t.A_1(x)^*,(\neg A_1(a))^*\ENTAIL.
\end{equation}
Using $a\leq t\ENTAIL Ea$ (a substitution instance of Axiom
\eqref{eq:Ax-E-p}), we obtain $E\vec{b},a\leq t,\forall x\leq
t.A_1(x)^*,(\neg A_1(a))^*\ENTAIL$ by Cut and Contraction. Using the
$\LEFTb\exists$-rule of $S^i_2E$, we obtain $E\vec{b},\forall x\leq
t.A_1(x)^*,\exists x\leq t.(\neg A_1(x))^*\ENTAIL$. By assumption,
$E\vec{a}, \Gamma^* \ENTAIL \Delta^*, \forall x\leq t.A_1(x)^*$ holds,
so by Cut and Contraction we have $E\vec{a},\exists x\leq t.(\neg
A_1(x))^*,\Gamma^*\ENTAIL\Delta^*$. By the definition of $*$, we are
done.

\paragraph{$A \equiv \forall x. A_1(x)$:} The proof is similar to the
proof of the previous case, the main differences being that we use
$Ea$ instead of $a\leq t$ and we apply the unbounded
counterparts of the bounded-quantifier rules of $S^i_2E$ employed in
that proof.

\paragraph{$A \equiv \exists x \leq t. A_1(x)$:} Let $a$ be a variable
that occurs in neither $A$ nor $\Gamma\ENTAIL\Delta$. As shown in the
proof of the case where $A\equiv\forall x\leq t.A_1(x)$, we can derive
$E\vec{b},Ea,A_1(a)^*,(\neg A_1(a))^*\ENTAIL$, where $\vec{b}$ are the
variables other than $a$ that occur free in $A_1(a)$. By the
$\LEFTb\forall$-rule of $S^i_2E$, we obtain $E\vec{b},Ea,a\leq t,
(A_1(a))^*,\forall x\leq t.(\neg A_1(x))^*\ENTAIL$.
Using $a\leq t\ENTAIL Ea$ (a substitution instance of Axiom
\eqref{eq:Ax-E-p}), we obtain $E\vec{b},a\leq t,(A_1(a))^*,\forall
x\leq t.(\neg A_1(x))^*\ENTAIL$ by Cut and Contraction. Using the
$\LEFTb\exists$-rule of $S^i_2E$, we obtain $E\vec{b}, \exists x\leq
t.A_1(x)^*,\forall x\leq t.(\neg A_1(x))^*\ENTAIL$. By assumption,
$E\vec{a},\Gamma^*\ENTAIL\Delta^*,\exists x\leq t.A_1(x)^*$ holds, so
by Cut and Contraction we have $E\vec{a},\forall x\leq t.(\neg
A_1(x))^*,\Gamma^*\ENTAIL\Delta^*$. By the definition of $*$, we are
done.

\paragraph{$A \equiv \exists x. A_1(x)$:} The proof is similar to the
proof of the previous case, the main differences being that we use
$Ea$ instead of $a \leq t$ and we apply the unbounded counterparts of
the bounded-quantifier rules of $S^i_2E$ used in that proof.

This completes the proof for the $\LEFT\neg$-rule.
\qed

\begin{lem}\label{lem:Rneg}
  \begin{equation}
    \infer{E\vec{a}, \Gamma^* \ENTAIL \Delta^*, (\neg A)^*}{
      E\vec{a}, A^*, \Gamma^* \ENTAIL \Delta^* }
  \end{equation} is admissible in $S^i_2E$, where $\vec{a}$ are the variables that occur
  free in the sequent $A,\Gamma\ENTAIL\Delta$.
\end{lem}

\proof
 By induction on $A$.
  
 If $A$ is atomic ($A\equiv p(t_1,t_2)$, where $p$ is $=$ or $\leq$),
 the inference given in the statement of the Lemma follows from the
 $\RIGHT\neg$-rule of $S^i_2E$, together with $E\vec{a} \ENTAIL Et_1$
 and $E\vec{a}\ENTAIL Et_2$, where the latter are derivable by
 Corollary \ref{cor:conv}.  In the induction step, we assume that
 $E\vec{a}, A^*, \Gamma^* \ENTAIL \Delta^*$ holds and show that
 $E\vec{a}, \Gamma^* \ENTAIL \Delta^*, (\neg A)^*$ holds. The proof
 depends on the form of $A$.

 \paragraph{$A \equiv \neg A_1$:} By Identity and Weakening, we obtain
 $E\vec{b}, (A_1)^* \ENTAIL (A_1)^*$, where $\vec{b}$ are the
 variables that occur free in $A_1$.  By the induction hypothesis
 applied to $A_1$, we have $E\vec{b} \ENTAIL (\neg A_1)^*, (A_1)^*$.
 By assumption, $E\vec{a}, (\neg A_1)^*, \Gamma^* \ENTAIL \Delta^*$ is
 derivable, so by Cut and Contraction we have $E\vec{a}, \Gamma^*
 \ENTAIL \Delta^*, (A_1)^*$.  By the definition of $*$, we are done
 (note that, by De Morgan duality, $(\neg\neg A_1)^*\equiv (A_1)^*$).

\paragraph{$A \equiv A_1 \AND A_2$:} Since $(A_1)^*, (A_2)^* \ENTAIL
(A_1)^* \AND (A_2)^*$ is derivable, together with the assumption that
$E\vec{a}, (A_1)^*\AND (A_2)^*, \Gamma^* \ENTAIL \Delta^*$ we obtain $E\vec{a}, (A_1)^*,
(A_2)^*, \Gamma^* \ENTAIL \Delta^*$.  By the
induction hypothesis applied twice in succession (once to $A_1$ and
once to $A_2$), we have $E\vec{a}, \Gamma^* \ENTAIL \Delta^*, (\neg
A_1)^*, (\neg A_2)^*$.  By the $\RIGHT\OR_1$- and $\RIGHT\OR_2$-rules
of $S^i_2E$, together with Contraction, $E\vec{a}, \Gamma^* \ENTAIL
\Delta^*, (\neg A_1)^* \OR (\neg A_2)^*$ is derivable.  By the
definition of $*$, we are done.
  
\paragraph{$A \equiv A_1 \OR A_2$:} Since $(A_1)^*\ENTAIL (A_1)^* \OR
(A_2)^*$ and $(A_2)^* \ENTAIL (A_1)^* \OR (A_2)^*$ are derivable, we
have $E\vec{a}, (A_1)^*, \Gamma^* \ENTAIL \Delta^*$ and $E\vec{a},
(A_2)^*, \Gamma^* \ENTAIL \Delta^*$ from the assumption that
$E\vec{a}, (A_1)^*\OR (A_2)^*, \Gamma^* \ENTAIL \Delta^*$.  By the
induction hypothesis applied to $A_1$ and $A_2$, we have $E\vec{a},
\Gamma^* \ENTAIL \Delta^*, (\neg A_1)^*$ and $E\vec{a}, \Gamma^*
\ENTAIL \Delta^*, (\neg A_2)^*$.  Thus we obtain $E\vec{a}, \Gamma^*
\ENTAIL \Delta^*, (\neg A_1)^* \AND (\neg A_2)^*$ by the
$\RIGHT\AND$-rule of $S^i_2E$.  By the definition of $*$, we are done.
  
\paragraph{$A \equiv A_1 \IMPLY A_2$:} Since $(\neg A_1)^* \ENTAIL
(\neg A_1)^* \OR (A_2)^*$ and $(A_2)^* \ENTAIL (\neg A_1)^* \OR
(A_2)^*$, we obtain $E\vec{a}, (\neg A_1)^*, \Gamma^* \ENTAIL
\Delta^*$ and $E\vec{a}, (A_2)^*, \Gamma^* \ENTAIL \Delta^*$ from the
assumption that $E\vec{a}, (\neg A_1)^*\OR (A_2)^*, \Gamma^* \ENTAIL
\Delta^*$.  By the induction hypothesis applied to $\neg A_1$ and
$A_2$, we have $E\vec{a}, \Gamma^* \ENTAIL \Delta^*, (A_1)^*$ (since
$(\neg \neg A_1)^* \equiv (A_1)^*$) and $E\vec{a}, \Gamma^* \ENTAIL
\Delta^*, (\neg A_2)^*$.  By the $\RIGHT\AND$-rule of $S^i_2E$, we
obtain $E\vec{a}, \Gamma^* \ENTAIL \Delta^*, (A_1)^* \AND (\neg
A_2)^*$.  By the definition of $*$, we are done.
\paragraph{$A \equiv \forall x \leq t. A_1(x)$:} Let $a$ be a variable
that occurs in neither $A$ nor $\Gamma\ENTAIL\Delta$. By Identity, we
have $A_1(a)^*\ENTAIL A_1(a)^*$. Thus $E\vec{b},Ea,a=a
,A_1(a)^*\ENTAIL A_1(a)^*$ holds by Weakening, where $\vec{b}$ are the
variables other than $a$ that occur free in $A_1(a)$. By the induction
hypothesis applied to $A_1(a)$, we have $E\vec{b},Ea,a=a\ENTAIL
A_1(a)^*, (\neg A_1(a))^*$. Using $Ea\ENTAIL a=a$ (Axiom
\eqref{eq:Ax-eq-refl}), we obtain $E\vec{b},Ea\ENTAIL A_1(a)^*,(\neg
A_1(a))^*$ by Cut and Contraction. By the $\RIGHTb\exists$-rule of
$S^i_2E$, we obtain $E\vec{b},Ea,a\leq t\ENTAIL A_1(a)^*,\exists x\leq
t.(\neg A_1(x))^*$.
Using $a\leq t\ENTAIL Ea$ (a substitution instance of Axiom
\eqref{eq:Ax-E-p}), we obtain $E\vec{b},a\leq t\ENTAIL
A_1(a)^*,\exists x\leq t.(\neg A_1(x))^*$ by Cut and
Contraction. Using the $\RIGHTb\forall$-rule of $S^i_2E$, we obtain
$E\vec{b},Et\ENTAIL \forall x\leq t.A_1(x)^*,\exists x\leq t.(\neg
A_1(x))^*$.  By assumption, $E\vec{a}, \forall x\leq t.(A_1(x))^*, \Gamma^* \ENTAIL
\Delta^*$ holds where $\vec{a}$ are free variables occuring in
$\forall x\leq t.(A_1(x)), \Gamma \ENTAIL \Delta$.  By Cut and
Contraction we have $E\vec{a},Et, \Gamma^* \ENTAIL \Delta^*, \exists
x\leq t. (\neg A_1(x))^*$. Using $E\vec{a}\ENTAIL Et$ (Corollary
\ref{cor:conv}), we obtain $E\vec{a},\Gamma^*\ENTAIL\Delta^*, \exists
x\leq t.(\neg A_1(x))^*$ by Cut and Contraction. By the definition of
$*$, we are done.

\paragraph{$A \equiv \forall x. A_1(x)$:} The proof is similar to the
proof of the previous case, the main differences being that we use
$Ea$ instead of $a \leq t$ and we apply the unbounded counterparts of
the bounded-quantifier rules of $S^i_2E$ used in that proof.
  
\paragraph{$A \equiv \exists x \leq t. A_1(x)$:} Let $a$ be a variable
that occurs in neither $A$ nor $\Gamma\ENTAIL\Delta$. As shown in the
proof of the case where $A\equiv \forall x\leq t.A_1(x)$, we can
derive $E\vec{b},Ea\ENTAIL A_1(a)^*,(\neg A_1(a))^*$, where $\vec{b}$
are the variables other than $a$ that occur free in $A_1(a)$. By the
$\RIGHTb\exists$-rule of $S^i_2E$, we obtain $E\vec{b},Ea,a\leq
t\ENTAIL \exists x\leq t.A_1(x)^*,(\neg A_1(a))^*$.  Using $a\leq
t\ENTAIL Ea$ (a substitution instance of Axiom \eqref{eq:Ax-E-p}), we
obtain $E\vec{b},a\leq t\ENTAIL \exists x\leq t.A_1(x)^*,(\neg
A_1(a))^*$ by Cut and Contraction. Using the $\RIGHTb\forall$-rule of
$S^i_2E$, we obtain $E\vec{b},Et\ENTAIL\exists x\leq
t.A_1(x)^*,\forall x\leq t.(\neg A_1(x))^*$. By assumption,
$E\vec{a},\exists x\leq t.(A_1(x))^*,\Gamma^*\ENTAIL\Delta^*$ holds,
so by Cut and Contraction we have
$E\vec{a},Et,\Gamma^*\ENTAIL\Delta^*,\forall x\leq t.(\neg
A_1(x))^*$. Using $E\vec{a}\ENTAIL Et$ (Corollary \ref{cor:conv}), we
obtain $E\vec{a},\Gamma^*\ENTAIL\Delta^*, \forall x\leq t.(\neg
A_1(x))^*$ by Cut and Contraction. By the definition of $*$, we are
done.

\paragraph{$A \equiv \exists x. A_1(x)$:} The proof is similar to the
proof of the previous case, the main differences being that we use
$Ea$ instead of $a\leq t$ and we apply the unbounded counterparts of
the bounded-quantifier rules used in that proof.
 
 This completes the proof for the $\RIGHT\neg$-rule.
\qed

\proof[Proof of Proposition \ref{prop:pred}]
  We prove that the $*$-translation of every inference of predicate
  logic is admissible in $S^i_2E$.  We consider only the rules of
  inference for negation, implication, and quantification.  The others
  are obvious.

  \paragraph{$\LEFT\neg$-rule:} This rule is treated in Lemma
  \ref{lem:Lneg}.
  
  \paragraph{$\RIGHT\neg$-rule:} This rule is treated in Lemma
  \ref{lem:Rneg}.

\paragraph{$\LEFT\IMPLY$-rule:}
\begin{equation}
  \infer{E\vec{a}, (\neg A_1)^* \OR A_2^*, \Gamma^* \ENTAIL \Delta^*
  }{
    E\vec{a}, \Gamma^* \ENTAIL \Delta^*, A_1^* \quad E\vec{a}, A_2^*,
    \Gamma^* \ENTAIL \Delta^*
  }
\end{equation} The admissibility of this inference follows from
Lemma \ref{lem:Lneg} and the $\LEFT\OR$-rule
of $S^i_2E$.

\paragraph{$\RIGHT\IMPLY$-rule:}
\begin{equation}
  \infer{E\vec{a}, \Gamma^* \ENTAIL \Delta^*,
    (\neg A_1)^* \OR (A_2)^*}{ E\vec{a}, (A_1)^*, \Gamma^* \ENTAIL
    \Delta^*, (A_2)^*
  }
\end{equation} The admissibility of this inference follows from Lemma
\ref{lem:Rneg} and the $\RIGHT\OR_1$- and
$\RIGHT\OR_2$-rules of $S^i_2E$, together with Contraction.

\paragraph{$\LEFTb\forall$-rule:}
\begin{equation}
  \infer[,]{E\vec{a}, E\vec{b}, t \leq s, \forall x \leq s. A(x)^*,
    \Gamma^* \ENTAIL \Delta^*}{ E\vec{a}, A(t)^*, \Gamma^* \ENTAIL
    \Delta^*
  }
\end{equation} where the variable $x$
does not occur in the term $s$, $\vec{a}$ are the variables that
occur free in the sequent $A(t),\Gamma\ENTAIL\Delta$, and $\vec{b}$
are the variables that occur in $s$ but are not in $\vec{a}$.
Without loss of generality, we can assume that $x$ does
not occur in $t$.  Therefore, $\vec{a}, \vec{b}$ are precisely
the variables that occur free in $t \leq s, \forall x
\leq s.A(x), \Gamma \ENTAIL \Delta$.

From the premise, we can derive $E\vec{a}, E\vec{b}, t \leq s,
\forall
x \leq s. A(x)^*, \Gamma^* \ENTAIL \Delta^*$ by the
$\LEFTb\forall$-rule of $S^i_2E$ and Weakening.

  \paragraph{$\LEFT\forall$-rule:}\begin{equation}
    \infer[,]{E\vec{a}, \forall x. A(x)^* , \Gamma^* \ENTAIL
    \Delta^*}{
      E\vec{a'}, A(t)^*, \Gamma^* \ENTAIL \Delta^* }
  \end{equation} where $\vec{a'}$ are the variables that occur free in
  the sequent $A(t),\Gamma\ENTAIL\Delta$ and $\vec{a}$ are the
  variables that occur free in the sequent $\forall
  x.A(x),\Gamma\ENTAIL\Delta$.

  Clearly, every variable that occurs free in $\forall x. A(x) ,
  \Gamma \ENTAIL \Delta$ also occurs free in $A(t), \Gamma \ENTAIL
  \Delta$. If there is at least one variable in $\vec{a'}$ which is
  not in $\vec{a}$, then for each such variable $b$, we substitute $0$
  for $b$ in both $Eb$ and $t$.  After repeated application of Cut
  with the axiom $\ENTAIL E0$, we obtain $E\vec{a}, A(t')^*, \Gamma^*
  \ENTAIL \Delta^*$, where $t'$ is obtained by substituting $0$ in $t$
  for every variable in $\vec{a'}$ which is not in $\vec{a}$.
  
  By the $\LEFT\forall$-rule of $S^i_2E$, we have $E\vec{a}, Et',
  \forall x. A(x)^*, \Gamma^* \ENTAIL \Delta^*$.  By Corollary
  \ref{cor:conv}, we have $E\vec{a} \ENTAIL Et'$.  Hence, we obtain
  $E\vec{a}, \forall x. A(x)^*, \Gamma^* \ENTAIL \Delta^*$ by Cut and
  Contraction.

\paragraph{$\RIGHTb\forall$-rule:}
\begin{equation}
  \infer[,]{E\vec{a}, \Gamma^* \ENTAIL \Delta^*, \forall x \leq
  t. A(x)^*}{
    E\vec{a}, Ea, a \leq t, \Gamma^* \ENTAIL \Delta^*, A(a)^* }
\end{equation}where neither the variable $a$ nor the
variable $x$ occurs in the term $t$; $a$ does not occur free in
$\Gamma\ENTAIL\Delta$; and $\vec{a}$ are the variables
other than $a$ that occur free in the sequent $a \leq t, \Gamma
\ENTAIL \Delta, A(a)$. Clearly, $\vec{a}$ are precisely the
variables that occur free in the sequent $\Gamma \ENTAIL \Delta,
\forall x \leq t. A(x)$.

Using $a \leq t \ENTAIL Ea$ (a substitution instance of Axiom
\eqref{eq:Ax-E-p}), we can eliminate the $Ea$ in the antecedent of the
premise by Cut and Contraction.  Therefore, we have $E\vec{a}, a \leq
t, \Gamma^* \ENTAIL \Delta^*, A(a)^*$.  Using the
$\RIGHTb\forall$-rule of $S^i_2E$, we can derive $E\vec{a}, Et,
\Gamma^* \ENTAIL \Delta^*, \forall x \leq t. A(x)^*$.  However, since
$Et$ is derivable from $E\vec{a}$, we obtain $E\vec{a}, \Gamma^*
\ENTAIL \Delta^*, \forall x \leq t. A(x)^*$.

\paragraph{$\RIGHT\forall$-rule:}
\begin{equation}
  \infer[,]{E\vec{a}, \Gamma^* \ENTAIL \Delta^*, \forall x. A(x)^*}{
    E\vec{a} \{, Ea\}, \Gamma^* \ENTAIL \Delta^*, A(a)^* }
\end{equation}where the variable $a$ does not occur free in
$\Gamma\ENTAIL\Delta$, and $\vec{a}$ are the variables other than $a$
that occur free in the sequent $\Gamma \ENTAIL \Delta, A(a)$. Clearly,
$\vec{a}$ are precisely the variables that occur free in the
sequent $\Gamma \ENTAIL \Delta, \forall x. A(x)$.
Here, the formula $Ea$ is enclosed in braces to indicate that it is
not included in the premise of the $*$-translation of the
$\RIGHT\forall$-rule unless the variable $a$ occurs free in $A(a)$.
Otherwise, $Ea$ can be added to the premise by % <- blank before ","
Weakening.  In either case, we can derive $E\vec{a},Ea,\Gamma^*
\ENTAIL \Delta^*, \forall x.A(x)^*$ by the $\RIGHT\forall$-rule of
$S^i_2E$.
  
\paragraph{$\exists$-rules:} The proofs of admissibility of the
$*$-translations of the $\exists$-rules are
analogous to the proofs of admissibility of the
$*$-translations of the $\forall$-rules.
\qed

\subsection{Bootstrapping Phase IV : $*$-translation of
  $\Sigma^b_i\text{-}\PIND$ rule is admissible in $S^i_2E$}
\label{subsec:bootstrap_IV}

Finally, we prove admissibility of the $*$-translation of the
$\Sigma^b_i\text{-}\PIND$ rule of $S^i_2$.

First, we prove that our formulation of $\PIND$, which uses the binary
successor functions, proves Buss' formulation of $\PIND$
\cite{Buss:Book}, which uses $\lfloor x / 2 \rfloor$.
\begin{lem}\label{lem:PIND}
  Assume that $\Gamma, Ea, A(\lfloor \frac{1}{2} a \rfloor) \ENTAIL
  A(a), \Delta$ is provable in $S^i_2E$, where the variable $a$ does
  not occur free in $\Gamma\ENTAIL\Delta$ and $A(a)$ is a $\Sigma^b_i$
  formula. Then $\Gamma, E\vec{a}, A(0) \ENTAIL A(t), \Delta$ is also
  provable in $S^i_2E$, where $\vec{a}$ are the variables that occur
  in the term $t$.
\end{lem}

\proof
  Note that $\lfloor \frac{1}{2} s_0a \rfloor = \lfloor \frac{1}{2}
  s_1a \rfloor = a$ if $Ea$ holds.  Therefore, substituting $s_0a$ and
  $s_1a$ for $a$ in $\Gamma, Ea, A(\lfloor \frac{1}{2} a \rfloor)
  \ENTAIL A(a), \Delta$ and applying Cut with $Ea \ENTAIL Es_0a$ and
  $Ea \ENTAIL Es_1a$, we obtain $\Gamma, Ea, A(a) \ENTAIL A(s_0a),
  \Delta$ and $\Gamma, Ea, A(a) \ENTAIL A(s_1a), \Delta$,
  respectively.  Combining $\Gamma, A(0) \ENTAIL A(0), \Delta$ and the
  $\Sigma^b_i\text{-}\PIND\text{-}E$ rule, we have $\Gamma, Et, A(0)
  \ENTAIL A(t), \Delta$.  Since $Et$ is derivable from $E\vec{a}$
  (Corollary \ref{cor:conv}), we have $\Gamma, E\vec{a}, A(0) \ENTAIL
  A(t), \Delta$.
\qed

\begin{prop}\label{prop:PIND}
  The $*$-translation of the $\PIND$ rule of $S^i_2$, that is,
  the inference
\begin{equation}\label{eq:PIND-*}
  \infer[,]{E\vec{a}\{,E\vec{b}\}, \Gamma^*, A(0)^* \ENTAIL A(t)^*,
  \Delta^*}{
    E\vec{a} \{,Ea\}, \Gamma^*, A(\lfloor a/2 \rfloor)^* \ENTAIL
    A(a)^*, \Delta^*
  }
\end{equation}is admissible in $S^i_2E$, where the variable $a$ does
not occur free in $\Gamma\ENTAIL\Delta$, $A(a)$ is a $\Sigma^b_i$
formula, $\vec{a}$ are the variables other than $a$ that occur free
in $\Gamma, A(\lfloor
a/2 \rfloor) \ENTAIL A(a), \Delta$, and $\vec{b}$ are the variables
that occur in $t$ but are not in $\vec{a}$.

The formula $Ea$ (in the antecedent of $E\vec{a} \{,Ea\}, \Gamma^*,
A(\lfloor a/2 \rfloor)^* \ENTAIL A(a)^*, \Delta^*$) is enclosed in
braces, as is $E\vec{b}$ (in the antecedent of $E\vec{a}
\{,E\vec{b}\}, \Gamma^*, A(0)^* \ENTAIL A(t)^*, \Delta^*$), to
indicate that $Ea$ and $E\vec{b}$ are not included in those
antecedents unless the variable $a$ occurs free in $A(a)$.
\end{prop}

\proof
  If $a$ does not occur free in $A(a)$, then the premise and the
  conclusion of \eqref{eq:PIND-*} are identical\, hence
  \eqref{eq:PIND-*} is admissible.

If $a$ occurs free in $A(a)$, then since $A(a)^*$ is a
  $\Sigma^b_i$ formula, we can apply Lemma \ref{lem:PIND} to obtain
  $E\vec{a},E\vec{b}, \Gamma^*, A(0)^* \ENTAIL A(t)^*,
  \Delta^*$.
\qed

Finally, we have the tools to prove Theorem \ref{thm:bootstrap}, which
was the main objective of this section.
\proof[Proof of Theorem \ref{thm:bootstrap}]
  By Propositions \ref{prop:equality} and \ref{prop:BASIC}, the
  $*$-translations of the axioms of $S^i_2$ are provable in $S^i_2E$.
  Furthermore, Propositions \ref{prop:pred} and \ref{prop:PIND}
  guarantee admissibility of the $*$-translations of the inferences of
  $S^i_2$.  Therefore, the $*$-translation of any sequent $\Gamma
  \ENTAIL \Delta$ which is provable in $S^i_2$ is provable in
  $S^i_2E$.

  In addition, the proofs of these propositions show that if a proof
  of an $S^i_2$ sequent $\Gamma \ENTAIL \Delta$ contains only
  $\Sigma^b_i$ and $\Pi^b_i$ formulae of $S^i_2$, then there is a
  proof of $(\Gamma \ENTAIL \Delta)^*$ that contains only $\Sigma^b_i$
  and $\Pi^b_i$ formulae of $S^i_2$E.
\qed

\section{Finitistic G\"{o}del Sentences of $S^{-1}_2E$}
\label{sec:goedel-sentence}

In this section, we investigate finitistic G\"{o}del sentences of
$S^{-1}_2E$.  Throughout this section, $i$ denotes a positive integer,
and $\lceil A(\underline{x})\rceil$ denotes the G\"odel number of the
formula obtained from $A(x_1)$ by substituting for the variable $x_1$
the numeral representation of the natural number $x$.  The purpose of
underscoring the $x$ in $\lceil A(\underline{x}) \rceil$ is to
indicate that the entity which is substituted is the numeral
\emph{representation} of the natural number $x$ (in the
meta-language), and not the \emph{variable} $x$ (in the object
language).

Now, let us motivate our investigation.  The most interesting question
concerning $S^{-1}_2 E$ is whether $S^{i}_2$ proves $i\Con$ or not.
If the answer is negative, we have $S^i_2 \not= S^{i+2}_2$, hence we
can conclude that the hierarchy $S^1_2,S^2_2,S^3_2,\dots$ does not
collapse.

Buss and Ignjatovi\'c \cite{BussUnprovability} used a finitistic
G\"odel sentence, together with Solovay's induction speed-up method,
to show that $S^i_2$ does not prove the consistency of proofs that are
comprised entirely of $\Sigma^b_i$ and $\Pi^b_i$ formulae and use only
BASIC axioms and the rules of inference of predicate logic.  Since
their notion of proofs and $i$-normal proofs have a certain
similarity, it looks as though it would be worthwhile to emulate their
method to prove $S^i_2 \not\vdash i\Con$.  Unfortunately, the induction
speed-up method in the form in which Buss and Ignjatovi\'c used it
does not work for $S^{-1}_2E$.  However, it would still be interesting
to investigate these G\"{o}del sentences.

First, let us see why the induction speed-up method does not work.
Buss and Ignjatovi\'c employed two approaches to the induction
speed-up method.  In their first approach, they bounded (provably in
$S^1_2$) the size of an $S^{-1}_2$ proof of a formula
$\phi(\underline{x})$ by $(x \# (x \# x))^m + n$, where $m$ and $n$
depend on the size of an $S^i_2$ proof of $\forall x. \phi(x)$. In
their other approach, which they applied to $\PV$ and $\PV^{-}$ (the
induction-free fragment of $\PV$), they presented a
polynomial-time computable function which converts a $\PV$ proof of a
sequent $\Gamma \ENTAIL \Delta$ with \emph{numerically restricted
  variables} (see p.\ 241 of \cite{BussUnprovability}) to a $\PV^{-}$
proof of that sequent.  We will use Proposition \ref{prop:soundness}
to show that neither of these approaches works for $S^{-1}_2E$.

Consider the first approach.  If it does work for
$S^i_2E$, then for every $i$-normal formula $\phi(x)$
that satisfies $S^i_2E \vdash \forall x. \phi(x)$, there
are natural numbers $m, n$ and $k$ that satisfy the
condition
\begin{equation}\label{eq:speedup}
  S^1_2 \vdash \forall x
  \exists w \leq t_k(x)^m+n. i\Prf(w, \lceil \phi(\underline{x})
  \rceil),
\end{equation}
where $t_k(x) \equiv \overbrace{x \# \ldots \# x}^k$.  Let $l > k$,
and derive a contradiction by considering strictly $i$-normal proofs
of $t_l(a) = t_l(a)$.  By \eqref{eq:speedup}, for every natural number
$d$ there exists a strictly $1$-normal (hence a strictly $i$-normal)
proof $w(d)$ of the formula $t_l(d) = t_l(d)$ such that $w(d) \leq
t_k(d)^m+n$.  By Proposition \ref{prop:soundness}, $T_i(w(d), \lceil
t_l(d) = t_l(d) \rceil, \rho)$ holds, where $\rho$ is the empty
environment.  By Lemma \ref{lem:T0-refl}, $w(d)> t_l(d)$.  However,
$t_l(d) > t_k(d)^m+n$ for sufficiently large $d$.  This contradicts
\eqref{eq:speedup}.

Now, consider the second approach.  Assume that it does work for
$S^1_2E$, that is, that there is a polynomial-time computable function
$f$ which converts $S^1_2E$ proofs of numerically restricted sequents
to strictly 1-normal proofs of those sequents.  Let $t_k(x)$ as in
\eqref{eq:speedup}, and consider the sequent $\ENTAIL Et_k(d)$ where
$d$ is a numeral.  Then $S^1_2E \vdash \ENTAIL Et_k(d)$, by
Proposition \ref{prop:total} and the fact that $S^1_2E \vdash Ed$.
The proof of $Et_k(d)$ can be taken to be of size $O(|d| \cdot 2^k)$.
Obviously, $\ENTAIL Et_k(d)$ is a sequent with numerically restricted
variables (since it has no variable).  Since induction speed-up works,
$f$ computes a strictly 1-normal (hence a strictly $i$-normal) proof
$w(d, k)$ of $\ENTAIL Et_k(d)$ in polynomial time.  Hence the size of
$w(d, k)$ is less than $O(|d|^c \cdot 2^{k \cdot c})$.  By Proposition
\ref{prop:soundness} $T(w(d, k), \lceil Et_k(d) \rceil, \rho)$ holds,
where $\rho$ is the empty environment.  By Lemma \ref{lem:T0-refl},
$w(d, k)\ge \rho(t(d))$.  Therefore, the size of $w(d, k)$ is greater
than $O(|d|^k)$.  However, for enough large $k$ and $d$, $O(|d|^c
\cdot 2^{k \cdot c}) < O(|d|^k)$. Contradiction.

Next, we investigate a countably infinite set of G\"odel sentences
$\forall x.\varphi_k(x)$ and the terms $t_k(x)\equiv \overbrace{x \#
  (x \# (\cdots (x \# x)))}^k$.  We show that $i\Con$ is equivalent to
$\forall x. \varphi_k(x)$ for every $k \geq 2$.  Therefore, we could
show that $S^i_2 \not\vdash i\Con$ by showing that $S^i_2 \not\vdash
\forall x.  \varphi_k(x)$ for some $k \geq 2$.

For $k\ge 1$, let $\varphi_k$ be a formula which
satisfies \begin{equation}
  \label{eq:phi} S^1_2 \vdash \forall x [ \varphi_k(x) \leftrightarrow
  \neg \exists w \leq t_k(x). i\Prf(w, \lceil \varphi_k(\underline{x})
  \rceil) ].
\end{equation}

We would like to prove that the G\"odel sentence $\forall
x.\varphi_k(x)$ is undecidable.  We can easily see that $S^i_2
\not\vdash \neg \forall x. \varphi_k(x)$.  For suppose that $S^i_2
\vdash \neg \forall x. \varphi_k(x)$ then $\forall x. \varphi_k(x)$ is
false by soundness of $S^i_2$, so there exists $d$ such that
$\varphi_k(d)$ is false.  By \eqref{eq:phi}, $\exists w \leq
t_k(d). i\Prf(w, \lceil \varphi_k(d) \rceil)$ is true.  Hence
$\varphi_k(d)$ has a strictly $i$-normal proof.  In particular, $S^i_2
\vdash \varphi_k(d)$. Since $\varphi_k(d)$ is false, this contradicts
the soundness of $S^i_2E$.

On the other hand, it looks as though it would be difficult to prove
that $S^i_2 \not\vdash \forall x. \varphi_k(x)$.  The crux of the
problem is that even if $S^i_2 \vdash \forall x. \varphi_k(x)$, there
is no (known) bound on the length of a strictly $i$-normal proof
$w(d)$ of $\varphi_k(d)$.  We know that $w(d)$ is not bounded by
$t_k(d)$, because by \eqref{eq:phi} that would contradict the
consistency of $S^i_2$.  However, perhaps $w(d)$ is bounded by
$t_l(d)$ for some $l>k$.

Finally, we prove that $S^i_2 \vdash i\Con \leftrightarrow \forall
x. \varphi_k(x)$ for $k \geq 2$.  The proof of $S^i_2 \vdash i\Con
\rightarrow \forall x. \varphi_k(x)$ uses a method similar to that
used in the proof of Theorem 4 on p.\ 135 of \cite{Buss:Book}.  The
proof of $S^i_2 \vdash i\Con \leftarrow \forall x. \varphi_k(x)$ is a
consequence of the trivial fact that a contradiction proves anything.

First, consider $S^i_2 \vdash i\Con \ENTAIL \forall x. \varphi_k(x)$.
By \eqref{eq:phi}, we have $S^1_2 \vdash \neg \varphi_k(x) \ENTAIL \exists w
\leq t_k(x). i\Prf(w, \lceil \varphi_k(\underline{x}) \rceil)$.  On
the other hand, using a method similar to that used for the proof of
Theorem 4 on p.\ 135 of \cite{Buss:Book}, we can prove that for every
$\Sigma^b_1$ formula $\psi(x)$, there is a term $u(x)$ such that
$S^1_2 \vdash \psi(x) \ENTAIL \exists w \leq u(x). 0\Prf(w, \lceil
\psi(\underline{x}) \rceil)$.  Since $\neg \varphi_k(x)$ is a
$\Sigma^b_1$ formula, there is a term $v(x)$ such that $S^1_2 \vdash
\neg \varphi_k(x) \ENTAIL \exists w \leq v(x). 0\Prf(w, \lceil \neg
\varphi_k(\underline{x}) \rceil)$.  Therefore, $S^1_2 \vdash \neg
\varphi_k(x) \ENTAIL \exists w. i\Prf(w, \lceil \ENTAIL \rceil)$.

Now consider $S^i_2 \vdash \forall x. \varphi_k(x) \ENTAIL i\Con$, and
assume that $i\Prf(w, \lceil \ENTAIL \rceil)$.  Then by Weakening
there exists $w'(x)=O(x)$ such that $i\Prf(w'(x), \lceil
\varphi_k(\underline{x}) \rceil)$.  Therefore, $w'(x) \leq t_k(x)$ for
sufficiently large $x$ if $k \geq 2$. (Note that for $k=1$, we have
$t_k(x)=x$, so this inequality does not hold.) Hence $\exists w' \leq
t_k(x). i\Prf(w', \lceil \varphi_k(\underline{x}) \rceil)$ for
sufficiently large $x$, in which case $S^i_2 \vdash \neg \varphi_k(x)$
by \eqref{eq:phi}. From this it follows that $S^i_2 \vdash \exists
x. \neg \varphi_k(x)$.

\section*{Acknowledgement}
  I thank Georgia Martin for thoroughly checking
  the drafts of this paper and proposing many improvements.
\bibliographystyle{plain}
\bibliography{my}

\end{document}